\documentclass{article}
\usepackage[utf8]{inputenc}
\usepackage{amsthm,amsmath,amssymb}
\usepackage{mathrsfs}
\usepackage{indentfirst}
\usepackage{natbib}
\usepackage{graphicx}
\usepackage{tikz-cd}
\usepackage{multicol}
\usepackage{multirow}
\usepackage{float}
\usepackage{geometry}
\usepackage{fancyhdr}
\usepackage{hyperref}
\usepackage{enumerate}
\usepackage{cite}
\usepackage{natbib}
\setcitestyle{numbers,square}

\title{\huge Ampleness of Automorphic Line Bundles on $U(2)$ Shimura Varieties}
\author{\Large Deding Yang}
\begin{document}

\newcounter{counter}
\setcounter{counter}{1}
\newtheorem{lemma}{Lemma}
\newtheorem{proposition}[lemma]{Proposition}
\newtheorem{theorem}[lemma]{Theorem}
\newtheorem{definition}[lemma]{Definition}
\newtheorem{remark}[lemma]{Remark}
\newtheorem{example}[lemma]{Example}
\newtheorem{assumption}[lemma]{Assumption}
\newtheorem{notation}[lemma]{Notation}

\maketitle

\medskip

\noindent\textbf{Abstract} Let $F$ be a totally real field in which $p$ is unramfied and let $S$ denote the integral model of the Hilbert modular variety with good reduction at $p$. Consider the usual automorphic line bundle $\mathcal{L}$ over $S$.  On the generic fiber, it is well known that $\mathcal{L}$ is ample if and only if all the coefficients are positive. On the special fiber, it is conjectured in \citep{Tian-Xiao} that $\mathcal{L}$ is ample if and only if the coefficients satisfy certain inequalities. We prove this conjecture for $U(2)$ Shimura varieties in this paper and deduce a similar statement for Hilbert modular varieties from this.

\tableofcontents

\section{Introduction}

\medskip
\noindent \textbf{Motivation}: 
Shimura varieties, which arise naturally as higher-dimensional analogues of modular curves, play a foundational role in the study of arithmetic algebraic geometry. The geometry of these varieties, which are moduli spaces parametrizing abelian varieties with additional structures is linked to multiple facets of number theory, automorphic forms and representation theory.

The study of the geometry of Shimura varieties in particularly useful in the Langlands program, establishing in-depth relations between the theory of Galois representations and the theory of automorphic representations, for example, the Eichler--Shimura relation. One way to make explicit this relation is to realize automorphic forms as cohomology of certain automorphic vector bundles on Shimura varieties. In particular, one may focus on the mod $p$ side and ask for analogous questions, such as the Serre weight conjecture and Artin conjecture for weight $1$ modular forms. One of the technical issues is to know whether certain cohomology group vanishes, or in other words, understanding the positivity of the corresponding automorphic vector bundles. 

The positivity of automorphic vector bundles is well known in characteristic $0$. For example, for Hilbert modular varieties, the modular line bundle $\omega^{\underline{t}}=\bigotimes\limits_{\tau\in\Sigma_{\infty}}\omega_\tau^{\otimes t_\tau}$ is ample if and only if $t_\tau>0$ for all $\tau$. Here $\Sigma_{\infty}$ denotes the set of infinite places of the totally real field $F$. We will revisit this fact later in the introduction part. It is then natural to ask for a criterion of ampleness in the characteristic $p$ case. This question was only partially answered prior to this paper (\citep{Tian-Xiao},\citep{Andreatta-Goren}). For example, Andreatta and Goren \citep{Andreatta-Goren} proved that for the special fiber $X$ of Hilbert modular surface (that is, the totally real field $F$ has degree $2$ over $\mathbb{Q}$, and $\Sigma_{\infty}=\{\tau_1,\tau_2\}$), the modular line bundle $\omega^{\underline{t}}$ is ample if and only if $pt_1>t_2>\frac{1}{p}t_1$ when $p$ is inert in $F$, and $t_1,t_2>0$ when $p$ splits in $F$. The aim of this paper is to fully solve the ampleness criterion of automorphic line bundles on the special fiber of Hilbert modular varieties and $U(2)$ Shimura varieties.

\medskip
\noindent\textbf{Setting:} We will primarily focus on $U(2)$ Shimura varieties in this paper. Let $F$ be a totally real field and $E/F$ be a CM quadratic extension in which $p$ is unramified, and each prime $\mathfrak{p}\subseteq \mathcal{O}_F$ above $p$ splits in $E$. Let $G$ be the $\mathbb{Q}$-algebraic group of simpletic similitudes associated to a two-dimensional $E$-vector space $V$ equipped with a nondegenerate alternating pairing $\langle\cdot,\cdot\rangle:V\times V\to \mathbb{Q}$ with $\langle\alpha x,y\rangle=\langle x,\bar\alpha y\rangle$ for all $\alpha\in E$ and $x,y\in V$. For $K=K^pK_p\subseteq G(\mathbb{A}_f)$ with $K^p$ sufficiently small and $K_p$ hyperspecial, we have an integral model $Sh_K(G)$ of the unitary Shimura variety. 

Let $X=Sh_K(G)_k$ denote the special fiber and $\mathcal{A}/X$ the universal abelian variety. $\mathcal{A}$ is equipped with an action of $\mathcal{O}_E$ and therefore we have a decomposition for the differential sheaf:
\begin{equation}
	\omega_{\mathcal{A}^\vee/X}=\bigoplus_{\tau}(\omega_{\mathcal{A}^\vee/X,\tilde\tau}\oplus \omega_{\mathcal{A}^\vee/X,\tilde\tau}).
\end{equation}
We assume, for simplicity, that the subbundles $\omega_{\mathcal{A}^\vee/X,\tilde\tau}$ and $\omega_{\mathcal{A}^\vee/X,\tilde\tau^c}$ are locally free of rank $1$ for all $\tau$ in the set $\Sigma_\infty$ of embeddings $F\hookrightarrow \mathbb{R}$ with $\tilde\tau,\tilde\tau^c$ the corresponding embeddings $E\hookrightarrow \mathbb{C}$ extending $\tau$. This would correspond to a special case of the signature condition. We will show in Section \ref{Unitary Shimura Varieties} that $[\omega_{\mathcal{A}^\vee/X,\tilde\tau}]=[\omega_{\mathcal{A}^\vee/X,\tilde\tau^c}]$ in $\textnormal{Pic}(X)_{\mathbb{Q}}$. We denote by $[\omega_\tau]$ to be the class of them. For any tuple $\underline{t}=(t_\tau)_{\tau\in\Sigma_{\infty}}\in \mathbb{Q}^{\Sigma_{\infty}}$, we define
\begin{equation}
	[\omega^{\underline{t}}]=\sum_{\tau\in\Sigma_{\infty}}t_\tau[\omega_\tau]\in\textnormal{Pic}(X)_\mathbb{Q}.
\end{equation}
This will be the core object of our study. We can now state the main theorem of this paper under the simplified signature condition above:
\begin{theorem}\label{Main Theorem}
	$[\omega^{\underline{t}}]$ is ample if and only if
	\begin{equation}
		pt_\tau>t_{\sigma^{-1}\tau}.
	\end{equation}
	for all $\tau\in \Sigma_{\infty}$. In particular, all $t_\tau$ are positive.
\end{theorem}

\medskip
\noindent
\textbf{Remark} (\emph{about the history})
\begin{enumerate}
	\item This ampleness criterion was conjectured by Tian--Xiao in \citep[\S ~6]{Tian-Xiao}, together with a proof of its necessity. In this paper, we prove the sufficiency. In addition, we believe that our proof is also valid for more general quaternionic Shimura varieties after some minor changes.
	\item For Hilbert modular surface, the ampleness criterion was proved by Andreatta and Goren \citep[\S~8]{Andreatta-Goren}. Our proof is a natural generalization of their result to higher dimensional cases using intersection theory and geometric tools.
	\item Brunebarbe, Goldring, Koskivirta and Stroh \citep{BGKS} announced a sufficient condition of ampleness of automorphic line bundles in a more general setting, but in the case of $U(2)$ Shimura varieties, their condition is not an equivalence condition.
\end{enumerate}

\medskip
\noindent
\textbf{Remark} (\emph{on the result})
\begin{enumerate}
	\item (Kodaira vanishing) Let $\underline{t}=(t_\tau)_{\tau\in\Sigma_{\infty}}$ be a tuple of positive integers. Let $d$ denote the dimension of $X$. The conditions in \citep[Corollary 2.8]{Deligne-Illusie} are satisfied for $p\gg0$, so we get
	\begin{equation}
		\begin{aligned}
			&H^j(X,\Omega_{X/k}^i\otimes \omega^{\underline{t}})=0 \  &\textnormal{for $i+j>d$},\\
			&H^j(X,\omega_{X/k}^i\otimes \omega^{-\underline{t}})=0 \  &\textnormal{for $i+j<d$}.
		\end{aligned}
	\end{equation}
	This allows us to apply our theorem (at least when $p$ is sufficiently large) to get vanishing of the cohomology of automorphic vector bundles.
	\item 	In Section 5, we relate the integral model of Hilbert modular varieties to that of unitary Shimura varieties. We deduce the ampleness criterion for modular line bundles over the special fiber $Y$ of the Hilbert modular variety directly from the above theorem:
	\begin{theorem}\label{Hilbert}
		Let $[\omega^{\underline{t}}]=\sum\limits_{\tau\in\Sigma_\infty}t_\tau[\omega_\tau]$ be the modular line bundle over $Y$. Then $[\omega^{\underline{t}}]$ is ample if and only if
		\begin{equation}
			pt_\tau> t_{\sigma^{-1}\tau}.
		\end{equation}
		for all $\tau\in \Sigma_{\infty}$. In particular, all $t_\tau$ are positive.	
	\end{theorem}
	\item Let $\mathcal{Y}$ be the canonical integral model over some $\mathcal{O}$ of the Hilbert modular variety with prime-to-$p$ level structure. Let $\underline{t}=(t_\tau)_{\tau\in\Sigma_{\infty}}$ be a tuple of positive numbers. It is known that $[\omega^{\underline{t}}]$ is ample over the generic fiber $\mathcal{Y}_\eta$. For this fixed tuple $\underline{t}$, if we choose $p$ sufficiently large, then the numerical condition in Theorem \ref{Hilbert} is automatically satisfied. We find that $[\omega^{\underline{t}}]$ is ample away from finitely many closed points on $\text{Spec } \mathcal{O}$. This corresponds to the gerenal fact in algebraic geometry that for a line bundle $\mathcal{L}$ on a relative scheme $Z$ over the base $S$, the set $\{s\in S\ |\ \mathcal{L}_s \text{ is ample on } Z_s\}$ is an open subscheme of $S$.
	\item We believe that our proof could be modified to give an ampleness criterion for splitting models of $U(2)$ Shimura varieties in the ramified case.
\end{enumerate}

\medskip
\noindent\textbf{Strategy of the proof}: We will use the following objects to deduce the ampleness criterion:
\begin{enumerate}
	\item (\emph{Goren--Oort stratification}) Goren and Oort \citep{Goren--Oort} defined a stratification on the special fiber of Hilbert modular varieties. Continue to assume that $\omega_{\mathcal{A}^\vee/X}$ and $\omega_{\mathcal{A}^\vee/X,\tilde\tau}$ are line bundles for all $\tau$ at this point. The Vershiebung map induces an $\mathcal{O}_E$-linear morphism $\omega_{\mathcal{A}^\vee/X}\longrightarrow\omega_{\mathcal{A}^{(p)\vee}/X}$, which further induces a homomorphism $h_{\tilde\tau}:\omega_{\mathcal{A}/X,\tilde\tau}\longrightarrow\omega_{\mathcal{A}^{(p)\vee}/X,\tilde\tau}\simeq \omega^{\otimes p}_{\mathcal{A}^\vee/X,\sigma^{-1}\tilde\tau}$, and similarly for the $\tilde\tau^c$-part for each $\tau\in\Sigma_{\infty}$. This map defines a global section $h_\tau\in H^0(X,\omega_{\mathcal{A}^\vee/X,\tilde\tau}^{\otimes -1}\otimes \omega_{\mathcal{A}^\vee/X,\sigma^{-1}\tilde\tau}^{\otimes p})$, which we call \emph{partial Hasse invariant at $\tilde\tau$}. We will prove in Section \ref{Zero loci} that $h_{\tilde\tau}$ and $h_{\tilde\tau^c}$ have the same scheme-theoretic zero loci. We denote by $X_\tau$ this common zero locus. For a subset $T\subset \Sigma_{\infty}$, we put $X_T=\cap_{\tau\in T}X_\tau$, and call it the \emph{Goren--Oort stratum of} $X$ \emph{at} $T$.
	
	\item (\emph{Description of the Goren--Oort strata}) Roughly speaking, each Goren--Oort stratum is isomorphic to a $\mathbb{P}^1$-power bundle over the special fiber of another $U(2)$ Shimura variety whose signature condition differs from the original one at some places. This perspective of relating different Shimura varieties can be traced back to Serre's two letters to Tate \citep{Serre-to-Tate} describing the supersingular locus of modular curves. This idea is widely used by many other authors in this area \citep{Diamond-Kassaei-Sasaki},\citep{Helm-Tian-Xiao},\citep{Tian-Xiao}. Our work relies mostly on Helm's isogeny trick in \citep{Helm}, where he studied sparse Goren--Oort strata. A complete description for all Goren--Oort strata is given in Tian-Xiao's work \citep{Tian-Xiao}. We slightly modify their proof and deduce the following theorem \textnormal{(Theorem \ref{description} in this paper)} in our case:
	\begin{theorem}
		The Goren--Oort stratum $X_T$ is isomorphic to a $(\mathbb{P}^1)^{I_T}$-bundle over $Sh_{K'}(G_{\mathcal{P}'})_k$. 
	\end{theorem}
	We will give the explicit definition of the terms in this theorem in Section \ref{Stratification}.
	
	\item (\emph{Inductive proof of the ampleness criterion}) We prove the necessity part of the ampleness criterion by restricting to the Goren--Oort stratum $X_{\Sigma_\infty\backslash\{\tau\}}$, which is a  $\mathbb{P}^1$-bundle over a discrete Shimura variety. Over each fiber, $[\omega^{\underline{t}}]=[\mathcal{O}(pt_\tau-t_{\sigma^{-1}\tau})]$, so we get the condition $pt_\tau>t_{\sigma^{-1}\tau}$.
	
	We prove the sufficiency part by induction. First, we write $[\omega^{\underline{t}}]=\epsilon[\det\omega]+[\omega^{\underline{t'}}]$ for a sufficiently small $\epsilon>0$. By the ampleness of the Hodge bundle $\det\omega$ \citep{Lan}, it suffices to prove the nefness of $[\omega^{\underline{t'}}]$ under the same numerical condition; Next, we rewrite $[\omega^{\underline{t}}]=\sum\limits_{\tau\in\Sigma_\infty}\lambda_\tau[h_\tau]$ and show that $\lambda_\tau\ge0$. So for a curve not contained in any strata, the intersection number $(C\cdot [\omega^{\underline{t}}])\ge0$; Finally, for $C$ contained in a stratum $X_T$ for some $T$, let $j:X_T\hookrightarrow X$ denote the natural embedding. By the description of Goren--Oort strata, there is proper morphism $\pi: X_T\longrightarrow Y$, where $Y$ is the special fiber of another $U(2)$-Shimura variety. Our strategy, roughly speaking, is to write $[\omega^{\underline{t}}]$ as a nonnegative linear combination of $\pi^\ast [\omega_Y^{\underline{t'}}]$ for some automorphic line bundle on $Y$ and some partial Hasse invariants. In the actual proof, we have to distinguish two cases ``the sparse case'' and ``the adjacent case'', and use several auxiliary Shimura varieties $Y_1,\dots, Y_k$. Then we can make induction and show the nonnegativity of the intersection $(C\cdot [\omega^{\underline{t}}])$ in this case.
	
\end{enumerate}

\medskip
\noindent
\textbf{Generalizations}: It is natural to ask for a generalization of this numerical result to other Hodge type Shimura varieties. To attack this problem, we want to understand the ampleness criterion for the automorphic line bundle $\mathcal{L}(\lambda)$ over the flag space $\pi:Fl\longrightarrow Sh_K(G)_k$ and relate it to the study of the automorphic vector bundle $\mathcal{V}(\lambda)$ on $Sh_K(G)_k$ using the Borel-Weil Theorem:
\begin{equation}
	\pi_\ast \mathcal{L}(\lambda)=\mathcal{V}(\lambda).
\end{equation}
We hope to address this question in some future works.

\medskip
\noindent
\textbf{Outline}: In Chapter \ref{Unitary Shimura Varieties}, we introduce basic facts of unitary Shimura varieties and partial Hasse invariants. Then we discuss the classes of certain automorphic line bundles in the rational Picard group, and give a concrete statement of the main theorem of this paper. In Chapter \ref{Stratification}, we define the Goren-Oort stratification and establish the geometric description of Goren--Oort strata, which realize the stratum $X_T$ as a $\mathbb{P}^1$-power bundle over the special fiber of another Shimura variety. A direct result is the necessity part of the ampleness criterion. Finally in Chapter \ref{Sufficiency}, we use induction and intersection theory to prove the sufficiency part of ampleness criterion.

\medskip
\noindent
\textbf{Acknowledgements}: First, I sincerely thank my advisor Liang Xiao for encouraging me to study the geometry of Shimura varieties. Under his guidance, I learned about relevant working methods in this field. Next, I want to thank Yichao Tian, Zhiyu Tian, Ruiqi Bai, Jiedong Jiang for many useful discussions. I would also like to express my thanks to them for encouraging me when I was stuck at some point.

\section{Unitary Shimura Varieties}\label{Unitary Shimura Varieties}

\subsection{Notations}

\begin{itemize}
	\item Let $p$ be a prime number. We fix an isomorphism $\iota: \mathbb{C}\stackrel{\sim}{\longrightarrow} \overline{\mathbb{Q}}_p$.
	\item Let $F$ be a totally real field of degree $d$ over $\mathbb{Q}$ in which $p$ is unramified. Let $\{\mathfrak{p}_1,\dots,\mathfrak{p}_g\}$ be the set of primes over $p$. For each $i$, define $f_i$ to be the inertia degree of $\mathfrak{p}_i$ over $p$. We denote by $\Sigma$ the set of places of $F$ and $\Sigma_\infty$ the subset of real embeddings. Via the isomorphism $\iota$, each $\tau\in \Sigma_\infty$ corresponds to a $p$-adic embedding $F\hookrightarrow \overline{\mathbb{Q}}_p$, and is therefore associated to a prime $\mathfrak{p}_i$ over $p$. We define $\Sigma_{\infty/\mathfrak{p}_i}$ to be the set of such embeddings. Note that $\Sigma_\infty=\bigsqcup_{i=1}^g \Sigma_{\infty/\mathfrak{p}_i}$.
	\item Since all $p$-adic embeddings take $F$ into the maximal unramified extension of $\mathbb{Q}_p$, for each $i$, we can label the embeddings in $\Sigma_{\infty/\mathfrak{p}_i}$ by $\tau_{i,1},\dots,\tau_{i,f_i}$ such that composing with the Frobenius endomorphism on $F_{\mathfrak{p}_i}$ changes $\tau_{i,j}$ to $\tau_{i,j+1}$.
	\item Let $E/F$ be a CM extension such that every $\mathfrak{p}_i$ splits in $E$, i.e., $\mathfrak{p}_i\mathcal{O}_E=\mathfrak{q}_i{\mathfrak{q}}_i^c$ for all $i$. We fix the choice of a prime $\mathfrak{q}_i$ over $\mathfrak{p}_i$ for each $i$. Define the set $\Sigma_E$ to be places of $E$ and its subsets $\Sigma_{E/\infty}$, $\Sigma_{E,\infty/\mathfrak{q}_i}$ and $\Sigma_{E,\infty/{\mathfrak{q}}_i^c}$ similarly as above. For each $\tau_i$, write $\tilde{\tau}_i$ for the embedding of $E$ extending $\tau_i$ that belongs to $\Sigma_{E,\infty/\mathfrak{q}_i}$, and $\tilde{\tau}^c_i$ the embedding of $E$ extending $\tau_i$ that belongs to $\Sigma_{E,\infty/{\mathfrak{q}}_i^c}$.
	\item Throughout this paper, we shall always make partitions $\mathcal{P}:\Sigma_\infty=\Sigma_{\infty,0}\bigsqcup\Sigma_{\infty,1}\bigsqcup\Sigma_{\infty,2}$ and $\Sigma_{\infty/\mathfrak{p},u}=\Sigma_{\infty,u}\cap \Sigma_{\infty/\mathfrak{p}}$ for $u=0,1,2$. Each partition will later correspond to a certain signature condition for some Shimura variety.
	\item Fix a finite field $k$ containing all residue fields of $F$ of characteristic $p$. Let $W(k)$ denote the ring of Witt vectors.
\end{itemize}

\subsection{Shimura datum}

\begin{itemize}
	\item Choose a two-dimensional $E$-vector space $V$ equipped with a nondegenerate alternating pairing $\langle\cdot,\cdot\rangle:V\times V\longrightarrow \mathbb{Q}$ satisfying
		\begin{equation}
				\langle \alpha x,y\rangle=\langle x,\bar\alpha y\rangle,\quad \alpha\in E.
		\end{equation}
and the following signature condition determined by the partition $\mathcal{P}$: For each $\tau\in \Sigma_{\infty}$, the pairing $\langle\cdot,\cdot\rangle$ induces an alternating pairing on the complex vector space $V_\tau=V\otimes_{F,\tau}\mathbb{R}$. This pairing is always the ``imaginary'' part of a uniquely determined Hermitian form $[\cdot,\cdot]_\tau$ on $V_\tau$. We require the signature of $[\cdot,\cdot]_\tau$ to be $(r_\tau,s_\tau)$ if $\tau\in\Sigma_{\infty,r_\tau}$ for some $r_\tau\in\{0,1,2\}$ and $s_\tau=2-r_\tau$.
	\item We define $G_\mathcal{P}$ to be the unitary similitude group of $V$. Concretely, for a $\mathbb{Q}$-algebra $R$,
		\begin{equation}
			G_{\mathcal{P}}(R)=\{(g,c(g))\in \mathrm{Aut}_{E\otimes_\mathbb{Q}R}(V\otimes_{\mathbb{Q}}R)\times R^\times,\  \langle gx,gy \rangle =c(g)\langle x,y\rangle,\  \forall x,y\in V\otimes_{\mathbb{Q}}R \}.
		\end{equation}
It is easy to see that
			\begin{equation}
				G_{\mathcal{P}}(\mathbb{R})=\mathrm{G}\bigg(\prod_{\tau\in \Sigma_\infty} \mathrm{U}(r_\tau,s_\tau)\bigg).
			\end{equation}
	\item We consider the homomorphism of $\mathbb{R}$-algebraic groups $h:\mathbb{S}=\rm{Res}_{\mathbb{C}/\mathbb{R}}\mathbb{G}_m\to G_{\mathcal{P},\mathbb{R}}$ defined by
			\begin{equation}
				h(z)=\prod_{\tau\in \Sigma_{\infty,0}} 
\left(\begin{matrix}
\bar{z} & 0\\
0 &\bar{z}
\end{matrix}\right)
\times \prod_{\tau\in \Sigma_{\infty,1}} 
\left(\begin{matrix}
z & 0\\
0 &\bar{z}
\end{matrix}\right)
\times \prod_{\tau\in \Sigma_{\infty,2}}
\left(\begin{matrix}
z & 0\\
0 &z
\end{matrix}\right).
			\end{equation}
for $z\in \mathbb{S}(\mathbb{R})=\mathbb{C}^\times$. Here the embedding $\textnormal{U}(r_\tau,s_\tau)\subseteq \textnormal{GL}_2(\mathbb{C})$ is induced by $\tilde\tau$. Now we have the isomorphism $\mathbb{S}_\mathbb{C}\simeq \mathbb{G}_{m,\mathbb{C}}\times\mathbb{G}_{m,\mathbb{C}}$, where the two factors are indexed by $\rm{Gal}(\mathbb{C}/\mathbb{R})\simeq\{\pm1\}$. Composing the complexified cocharacter $h_{\mathbb{C}}$ with the embedding of $\mathbb{G}_{m,\mathbb{C}}$ into the first factor gives a cocharacter $\mu_h:\mathbb{G}_{m,\mathbb{C}}\to G_{\mathcal{P},\mathbb{C}}$. Let $X$ be the $G_{\mathcal{P}}(\mathbb{R})$-conjugacy class of homomorphisms $\mu_h$.
	\item The above defines the Shimura datum $(G_\mathcal{P},X)$.

\end{itemize}

\subsection{Moduli interpretation}
	Let $K=K^pK_p\subseteq G(\mathbb{A}_\mathbb{Q}^f)$ be an open subgroup with hyperspecial level structure at $p$. In other words, $G$ extends to a reductive group scheme $\mathcal{G}/\mathbb{Z}_p$ with $\mathcal{G}\otimes_{\mathbb{Z}_p}\mathbb{Q}_p\simeq G_{\mathbb{Q}_p}$ and $K_p=\mathcal{G}(\mathbb{Z}_p)$. According to \citep{Kottwitz}
for $K^p$ sufficiently small, the Shimura variety $Sh_K(G_{\mathcal{P}})$ can be realized as the moduli scheme of the functor:
\begin{equation}
Sch^{loc\  noe}_{/W(k)}\longrightarrow Sets
\end{equation}
sending a locally noetherian scheme $S$ over $W(k)$ to the set $Sh_K(G_\mathcal{P})(S)$ consisting of equivalent classes of tuples $(A,\lambda,\eta)$ satisfying
\begin{enumerate}
	\item $A$ is an abelian scheme of dimension $2d$ over $S$ equipped with an action of $\mathcal{O}_E$, such that the characteristic polynomial of $\alpha\in\mathcal{O}_E$ on $\textnormal{Lie}_{A/S}$ is gived by
		\begin{equation}
			\prod_{\tau:F\to\mathbb{R}}(x-\tilde{\tau}(\alpha))^{r_\tau}(x-\tilde{\tau}^c(\alpha))^{s_\tau}
		\end{equation}
In other words, if $\textnormal{Lie}_{A/S,\tilde\tau}$ (resp. $\textnormal{Lie}_{A/S,\tilde\tau^c}$) denotes the subsheaf of $\textnormal{Lie}_{A/S}$ on which $\mathcal{O}_E$ acts via $\tilde\tau$ (resp. via $\tilde\tau^c$), then it is a subbundle of rank $r_\tau$ (resp. $s_\tau$).
	\item $\lambda:A\to A^\vee$ is a prime-to-$p$ quasi-polarization of $A$ such that the Rosati involution associated to $\lambda$ induces the complex conjugation on $\mathcal{O}_E$.
	\item $\eta$ is a $K^p$-level structure on $A$. Explicitly, this means that $\eta$ is a collection of $\pi_1(S_j,\bar{s}_j)$-invariant $K^p$-orbit of isomorphisms $\eta_j:V\otimes_\mathbb{Q}\mathbb{A}^{(p)}_f\stackrel{\sim}{\longrightarrow}V^{(p)}(A_{\bar{s}_j})$ for each connected component $S_j$ of $S$ with a geometric point $\bar{s}_j$. Here $V^{(p)}(A_{\bar{s}_j})= \varprojlim\limits_{p\nmid N}(A_{\bar{s}_j}[N])\otimes_\mathbb{Z}\mathbb{Q}$ is the rational Tate module away from $p$. We further require that for some isomorphism $\nu(\eta_j)\in \text{Hom}({\mathbb{A}}_f^{(p)},{\mathbb{A}}_f^{(p)}(1))$, the following diagram is commutative:

\centering{
\begin{tikzcd}
V\otimes_{\mathbb{Q}} \mathbb{A}_f^{(p)}\times V\otimes_{\mathbb{Q}} \mathbb{A}_f^{(p)} \arrow[rr, "{\langle\cdot,\cdot\rangle}"] \arrow[d, "\eta_j", shift left=10] \arrow[d, "\eta_j"', shift right=10] &  & {\mathbb{A}}_f^{(p)} \arrow[d, "\nu(\eta_j)"] \\
V^{(p)}(A_{\bar{s}_j})\times V^{(p)}(A_{\bar{s}_j}) \arrow[rr, "\textnormal{Weil pairing}"]                                                                                    &  & {\mathbb{A}}_f^{(p)}(1).              
\end{tikzcd}}
\end{enumerate}
Two triples $(A,\lambda,\eta)$ and $(A',\lambda',\eta')$ are equivalent if there is an $\mathcal{O}_E$-equivariant prime-to-$p$ quasi-isogeny $\phi:A\to A'$ such that $\phi^\vee\circ\lambda'\circ\phi=\lambda$ and $\eta'=\phi\circ \eta$.

For the rest of the paper, we denote by $X=Sh_K(G_\mathcal{P})_k$ the special fiber of the Shimura variety.

\subsection{Additional structures on $X$}

In this subsection, let $S$ be a noetherian $W(k)$-scheme and $(A,\lambda,\eta)$ be an $S$-point of $Sh_K(G_{\mathcal{P}})$.

\begin{itemize}
	\item (Hodge filtration) The first de Rham homology of an abelian scheme $A/S$ admits a filtration
		\begin{equation}\label{Hodge}
			0\longrightarrow \omega_{A^\vee/ S}\longrightarrow \mathrm{H}_1^{\mathrm{dR}}(A/S)\longrightarrow \text{Lie}_{A/S}\longrightarrow 0.
		\end{equation}
Thanks to the $\mathcal{O}_E$-action, we have a decomposition 
\begin{equation}
\textnormal{H}_1^{\textnormal{dR}}(A/S)=\bigoplus_{\tau\in\Sigma_\infty}(\textnormal{H}_1^{\textnormal{dR}}(A/S)_{\tilde\tau}\oplus \textnormal{H}_1^{\textnormal{dR}}(A/S)_{\tilde\tau^c}). 
\end{equation}
and similarly for $\omega_{A^\vee/S}$ and $\textnormal{Lie}_{A/S}$. Here the subscripts $\tilde\tau$ or $\tilde\tau^c$ indicate that $\mathcal{O}_E$ acts on the corresponding sheaf via $\tilde\tau$ or $\tilde\tau^c$. The exact sequence (\ref{Hodge}) decomposes into short exact sequences
\begin{equation}
0\longrightarrow \omega_{A^\vee/ S,\tilde\tau}\longrightarrow \mathrm{H}_1^{\mathrm{dR}}(A/S)_{\tilde\tau}\longrightarrow \text{Lie}_{A/S,\tilde\tau}\longrightarrow 0.
\end{equation}
and similarly for the $\tilde\tau^c$-components.
	\item (Frobenius and Verschiebung) Assume $S$ is a scheme of characteristic $p$, let $F=F_{A/S}$ be the relative Frobenius morphism. It induces a homomorphism 
		\begin{equation}
			V=V_\tau:\mathrm{H}^{\mathrm{dR}}_1(A/S)_{\tilde\tau}\longrightarrow \mathrm{H^{dR}_1}(A^{(p)}/S)_{\tilde\tau}\stackrel{\sim}{\longrightarrow} \mathrm{H}^{\mathrm{dR}}_1(A/S)_{\sigma^{-1}\tilde\tau}^{(p)}
		\end{equation}
for each $\tilde\tau$, and similarly for $\tilde\tau^c$. Here the symbol $^{(p)}$ means the pullback via the absolute Frobenius on $S$. Similarly, we have the Verscheibung morphism $A^{(p)}\to A$ such that $V\circ F=[p]_A$, and it induces the homomorphism
		\begin{equation}
			F=F_\tau:\text{H}^{\text{dR}}_1(A^{(p)}/S)_{\tilde\tau}\stackrel{\sim}{\longrightarrow}\text{H}^{\text{dR}}_1(A/S)^{(p)}_{\sigma^{-1}\tilde\tau}\longrightarrow \text{H}^{\text{dR}}_1(A/S)_{\tilde\tau}.
		\end{equation}
The following identifications will be important for us later
		\begin{equation}
			\text{Ker}\ V=\text{Im}\ F;\quad \text{Ker}\ F=\text{Im}\ V=\omega_{A^\vee/S,\sigma^{-1}\tilde\tau}^{(p)}.
		\end{equation}

	\item (Dieudonn$\acute{\text{e}}$ theory) Assume $S$\ =\ Spec $\ell$ with $\ell$ a perfect field containing $k$, let $W(\ell)$ denote the ring of Witt vectors of $\ell$. Let $\tilde{\mathbb{D}}(A)$ be the covariant Dieudonn$\acute{\text{e}}$ module associated to the $p$-divisible group $A[p^\infty]$ of $A$. This is a finite free $W(\ell)$-module equipped with a $\sigma$-linear action $F$ and a $\sigma^{-1}$-linear action $V$ satisfying $FV=VF=p$. The $\mathcal{O}_E$-action on $A$ naturally induces an $\mathcal{O}_E$-action on $\tilde{\mathbb{D}}(A)$ which commutes with the operators $F$ and $V$. Let $\tilde{\mathbb{D}}(A)_{\tilde\tau}$ and $\tilde{\mathbb{D}}(A)_{\tilde\tau^c}$ be the submodules of $\tilde{\mathbb{D}}(A)$ over which $\mathcal{O}_E$ acts via $\tilde\tau$ and $\tilde\tau^c$ respectively. We have the comparison
		\begin{equation}
			\tilde{\mathbb{D}}(A)/p\tilde{\mathbb{D}}(A)\stackrel{\sim}{\longrightarrow} \mathrm{H}^{\mathrm{dR}}_1(A/S) \quad \textnormal{and}\quad \tilde{\mathbb{D}}(A)_{\tilde\tau}/p\tilde{\mathbb{D}}(A)_{\tilde\tau}\stackrel{\sim}{\longrightarrow} \mathrm{H}^{\mathrm{dR}}_1(A/S)_{\tilde\tau}.
		\end{equation}
for each $\tilde\tau$ and similarly for each $\tilde{\tau}^c$.

	\item (Partial Hasse invariants) We continue to assume that $S$ is a scheme in characteristic $p$. Classically, Hasse invariants measure the behavior of the differentials under the Frobenius morphism. We consider partial Hasse invariants, which are ``factors'' of the classical Hasse invariants. For a place $\tau\in \Sigma_{\infty}$, we first introduce essential Frobenius and Verchiebung
\begin{equation}
F_{es}:\textnormal{H}_1^{\textnormal{dR}}(A^{(p)}/S)\to \textnormal{H}_1^{\textnormal{dR}}(A/S)\ \textnormal{and}\ V_{es}:\textnormal{H}_1^{\textnormal{dR}}(A/S)\to \textnormal{H}_1^{\textnormal{dR}}(A^{(p)}/S)
\end{equation}
defined on each $\tilde\tau$-component by
\begin{equation}
	F_{es,\tilde\tau}=\left\{
\begin{aligned}
F\quad  ,\ &\text{if } \sigma^{-1}\tau\in \Sigma_{\infty,1}\bigsqcup \Sigma_{\infty,2}  \\
V^{-1},\ &\text{if } \sigma^{-1}\tau\in \Sigma_{\infty,0}
\end{aligned}\right.
\ \text{and } V_{es,\tilde\tau}=\left\{
\begin{aligned}
V\quad,\ &\text{if }  \sigma^{-1}\tau\in \Sigma_{\infty,0}\bigsqcup\Sigma_{\infty,1}\\
F^{-1},\ &\text{if }  \sigma^{-1}\tau\in \Sigma_{\infty,2}
\end{aligned}\right..
\end{equation}
The upshot is that, both $F_{es,\tilde\tau}$ and $V_{es,\tilde\tau}$ are isomorphisms or have cokernels locally free of rank $1$, where the latter occurs if and only if $\sigma^{-1}\tilde\tau\in \Sigma_{\infty,1}$. Now if we start with a place $\tau\in\Sigma_{\infty,1}$, and let $n_\tau$ be the positive integer such that $\sigma^{-1}\tau,\ldots,\sigma^{-n_\tau+1}\tau\notin \Sigma_{\infty,1}$ but $\sigma^{-n_\tau}\tau\in \Sigma_{\infty,1}$, we get the composite
		\begin{equation}
			V_{es,\tilde\tau}^{n_\tau}:\text{H}^{\text{dR}}_1(A/S)_{\tilde\tau}\xrightarrow{V_{es,\tilde\tau}}\text{H}^{\text{dR}}_1(A/S)^{(p)}_{\sigma^{-1}\tilde\tau}\xrightarrow{V^{(p)}_{es,\sigma^{-1}\tilde\tau}}\cdots\xrightarrow{V_{es,\sigma^{1-n_\tau}\tilde\tau}^{(p^{n_\tau-1})}}\text{H}^{\text{dR}}_1(A/S)^{(p^{n_\tau})}_{\sigma^{-n_\tau}\tilde\tau},
		\end{equation}
where all the morphisms in the definition of $V^{n_\tau}_{es,\tilde\tau}$ are isomorphisms except the last one. Restricting this map to $\omega_{A^\vee/S,\tilde\tau}$ gives:
		\begin{equation}
			h_{\tilde\tau}(A):\omega_{A^\vee/S,\tilde\tau}\longrightarrow\omega_{A^{\vee,(p^{n_\tau})}/S,\tilde\tau}\stackrel{\sim}{\longrightarrow}\omega^{\otimes p^{n_\tau}}_{A^\vee/S,\sigma^{-n_\tau}\tilde\tau}.
		\end{equation}
Applying this to the universal abelian variety $\mathcal{A}$ over the special fiber $X=Sh_K(G_{\mathcal{P}})_k$, we get a section of a line bundle
		\begin{equation}\label{Hasse}
			h_{\tilde\tau}\in H^0(X,\omega_{\mathcal{A}^\vee/X,\tilde\tau}^{-1}\otimes\omega^{\otimes p^{n_\tau}}_{\mathcal{A}^\vee/X,\sigma^{-n_\tau}\tilde\tau}).
		\end{equation}
Similar definitions work for $\tilde\tau^c$ with all $0$ and $2$'s interchanged according to the signature conditions. We call these $h_{\tilde\tau},h_{\tilde\tau^c}$'s the $\tilde\tau,\tilde\tau^c$-\emph{partial Hasse invariants}, respectively.

\end{itemize}

\subsection{Zero loci of partial Hasse invariants}\label{Zero loci}

We compare the zero loci of $h_{\tilde\tau}$ and $h_{\tilde\tau^c}$ on the special fiber $X$ (see \citep{Tian-Xiao}, Lemma 4.5).

\begin{lemma} Let $(A,\lambda,\eta)\in X(S)$ be a closed point and $\tau\in \Sigma_{\infty,1}$. Then the following statements are equivalent:

	(1) $h_{\tilde\tau}(A)=0$.
	
	(2) \textnormal{Im}$\Big(F^{n_\tau}_{es,\tilde\tau}:\textnormal{H}^{\textnormal{dR}}_1(A^{(p^{n_\tau})}/S)_{\tilde\tau}\longrightarrow \textnormal{H}^{\textnormal{dR}}_1(A/S)_{\tilde\tau}\Big)=\omega_{A^\vee/S,\tilde\tau}$.

	(3) $h_{\tilde\tau^c}(A)=0$.
	
	(4) $\textnormal{Im}\Big(F^{n_\tau}_{es,\tilde\tau^c}:\textnormal{H}^{\textnormal{dR}}_1(A^{(p^{n_\tau})}/S)_{\tilde\tau^c}\longrightarrow \textnormal{H}^{\textnormal{dR}}_1(A/S)_{\tilde\tau^c}\Big)=\omega_{A^\vee/S,\tilde\tau^c}$.
\begin{proof} The fact that $\text{Im }F=\text{Ker }V$ implies the equivalences $(1)\Leftrightarrow(2)$ and $(3)\Leftrightarrow(4)$. To show $(2)\Leftrightarrow(4)$, recall that the Frobenius and Verschiebung maps on the de Rham cohomology groups are compatible with the (perfect) Weil pairing, in the sense that $\langle F^{n_\tau}_{es,\tilde\tau}x, y\rangle=\langle x, V^{n_\tau}_{es,\tilde\tau^c}y\rangle^{\sigma^{n_\tau}}$ for $x\in \textnormal{H}_1^{\textnormal{dR}}(A^{(p^{n_\tau})}/S)_{\tilde\tau}$ and $y\in \textnormal{H}^{\textnormal{dR}}_1(A/S)_{\tilde\tau^c}$. We therefore get the relations
	\begin{equation}
		\omega_{A^\vee/S,\tilde\tau}^\perp=\omega_{A^\vee/S,\tilde\tau^c}\quad \text{and} \quad \text{Im}(F^{n_\tau}_{es,\tilde\tau})^\perp=\text{Im}(F^{n_\tau}_{es,\tilde\tau^c}).
	\end{equation}
The equivalence of $(2)$ and $(4)$ then follows.
\end{proof}
\end{lemma}

\begin{notation}
We denote by $X_\tau$ to be the scheme-theoretic vanishing locus of $h_{\tilde\tau}$ (or $h_{\tilde\tau^c}$) on $X$.
\end{notation}

\begin{lemma}\label{Picard class} Let $\mathcal{A}/X$ be the universal abelian scheme. Then in the rational Picard group $\textnormal{Pic}(X)_{\mathbb{Q}}=\textnormal{Pic}(X)\otimes_{\mathbb{Z}}\mathbb{Q}$ of $X$, we have
	
	(1) If $\tau\in\Sigma_{\infty,1}$, then $[\omega_{\mathcal{A}^\vee/X,\tilde\tau}]=[\omega_{\mathcal{A}^\vee/X,\tilde\tau^c}]=[\omega_{\mathcal{A}/X,\tilde\tau}]=[\omega_{\mathcal{A}/X,\tilde\tau^c}]$.

	(2) For any $\tau\in\Sigma_{\infty}$, $[\wedge^2 \textnormal{H}^{\textnormal{dR}}_1(\mathcal{A}/X)_{\tilde\tau}]=[\wedge^2 \textnormal{H}^{\textnormal{dR}}_1(\mathcal{A}/X)_{\tilde\tau^c}]=0$.
\end{lemma}

	(3) Let $X^{\textnormal{min}}$ be the minimal compactification of $X$. Then the natural homomorphism $j:\textnormal{Pic}(X^{\textnormal{min}})\to\textnormal{Pic}(X)$ is injective. Moreover, both $[\omega_{\mathcal{A}^\vee/X,\tilde\tau}]$ and $[\omega_{\mathcal{A}^\vee/X,\tilde\tau^c}]$ belong in the image of $j_{\mathbb{Q}}:\textnormal{Pic}(X^{\textnormal{min}})_{\mathbb{Q}}\to \textnormal{Pic}(\textnormal{X})_{\mathbb{Q}}$.

\begin{proof} This is also an adaptation of Lemma 6.2 in \citep{Tian-Xiao}. First, the polarization on $\mathcal{A}$ induces isomorphisms $\omega_{\mathcal{A}^\vee/X,\tilde\tau}\stackrel{\simeq}{\longrightarrow}\omega_{\mathcal{A}/X,\tilde\tau^c}$ and $\omega_{\mathcal{A}/X,\tilde\tau}\stackrel{\simeq}{\longrightarrow}\omega_{\mathcal{A}^\vee/X,\tilde\tau^c}$. Next, by the previous lemma, the zero loci of $h_{\tilde\tau}$ and $h_{\tilde\tau^c}$ have the same underlying sets. We will prove in Section 3 that the zero loci of partial Hasse invariants are smooth, so $h_{\tau}$ and $h_{\tilde\tau}$ define the same Weil divisors. We have 
	\begin{equation}
		p^{n_\tau}[\omega_{\mathcal{A}^\vee/X,\sigma^{-n_\tau}\tilde\tau}]-[\omega_{\mathcal{A}^\vee/X,\tilde\tau}]=[h_{\tilde\tau}]=[h_{\tilde\tau^c}]=p^{n_\tau}[\omega_{\mathcal{A}^\vee/X,\sigma^{-n_\tau}\tilde\tau^c}]-[\omega_{\mathcal{A}^\vee/X,\tilde\tau^c}].
	\end{equation}
We may then define a matrix $H$ labelled by $\tau\in\Sigma_{\infty,1}$, with $(\tau_1,\tau_2)$-th entry given by
	\begin{equation}
		h_{\tau_1,\tau_2}=\left\{
\begin{aligned}
-1\ ,\quad &\text{if }\tau_1=\tau_2, \\
p^{n_{\tau_2}},\quad &\text{if }\tau_1=\sigma^{-n_{\tau_2}}\tau_2,\\
0\quad,\quad &\text{otherwise}.
\end{aligned}\right.
	\end{equation}
This matrix is the transformation matrix from the basis given by $\omega_{A^\vee/X,\tilde\tau}$ (or $\omega_{A^\vee/X,\tilde\tau^c}$) to the basis given by partial Hasse invariants. It is easy to see that $H$ is invertible, so the above equation implies $[\omega_{\mathcal{A}^\vee/X,\tilde\tau}]=[\omega_{\mathcal{A}^\vee/X,\tilde\tau^c}]$.

For the second statement, if $\tau\in \Sigma_{\infty,1}$, then both sides of the Hodge filtration 
\begin{equation}
	0\longrightarrow \omega_{\mathcal{A}^\vee/X,\tilde\tau}\longrightarrow \textnormal{H}_1^{\textnormal{dR}}(\mathcal{A}/X)_{\tilde\tau}\longrightarrow \text{Lie}_{\mathcal{A}/X,\tilde\tau}\longrightarrow0.
\end{equation}
are line bundles, so we have $[\wedge^2 \text{H}^{\text{dR}}_1(\mathcal{A}/X)_{\tilde\tau}]=[\omega_{\mathcal{A}^\vee/X,\tilde\tau}]+[\text{Lie}_{\mathcal{A}/X,\tilde\tau}]=[\omega_{\mathcal{A}^\vee/X,\tilde\tau}]-[\omega_{\mathcal{A}/X,\tilde\tau}]=0$; if $\tau\notin \Sigma_{\infty,1}$, first suppose that there is an integer $m$ such that $\sigma^{m}\tau\in\Sigma_{\infty,1}$ and $\sigma^{i}\tau\notin\Sigma_{\infty,1}$ for all $0\le i\le m-1$. We have a sequence of isomorphisms
\begin{equation}
\text{H}_1^{\text{dR}}(\mathcal{A}/X)_{\tilde\tau}^{(p^m)}\xrightarrow[\simeq]{F_{\mathcal{A},es}}\text{H}_1^{\text{dR}}(\mathcal{A}/X)_{\sigma\tilde\tau}^{(p^{m-1})}\xrightarrow[\simeq]{F_{\mathcal{A},es}}\cdots\xrightarrow[\simeq]{F_{\mathcal{A},es}}\text{H}^{\text{dR}}_1(\mathcal{A}/X)_{\sigma^m\tilde\tau}.
\end{equation}
This gives
\begin{equation}
p^{m}[\wedge^2 \text{H}^{\text{dR}}_1(\mathcal{A}/X)_{\tilde\tau}]=[\wedge^2 \text{H}^{\text{dR}}_1(\mathcal{A}/X)_{\sigma^m\tilde\tau}]=0.
\end{equation}
If such $m$ does not exist, then there exists $d$ such that $\sigma^{2d}\tilde\tau=\tilde\tau$. The $2d$-time composition of isomorphisms
\begin{equation}
\text{H}_1^{\text{dR}}(\mathcal{A}/X)_{\tilde\tau}^{(p^{2d})}\xrightarrow[\simeq]{F_{\mathcal{A},es}}\text{H}_1^{\text{dR}}(\mathcal{A}/X)_{\sigma\tilde\tau}^{(p^{2d-1})}\xrightarrow[\simeq]{F_{\mathcal{A},es}}\cdots\xrightarrow[\simeq]{F_{\mathcal{A},es}}\text{H}^{\text{dR}}_1(\mathcal{A}/X)_{\sigma^{2d}\tilde\tau}.
\end{equation}
gives the equation
\begin{equation}
p^{2d}[\wedge^2 \text{H}^{\text{dR}}_1(\mathcal{A}/X)_{\tilde\tau}]=[\wedge^2 \text{H}^{\text{dR}}_1(\mathcal{A}/X)_{\tilde\tau}].
\end{equation}
so $[\wedge^2 \text{H}^{\text{dR}}_1(\mathcal{A}/X)_{\tilde\tau}]=0$. Similar argument applies to $\tilde\tau^c$.

For the last statement, note that $X^{\textnormal{min}}$ is normal and $X^{\textnormal{min}}-X$ has codimension$\ge2$ by \citep{Lan}, so $j:\textnormal{Pic}(X^{\textnormal{min}})\to \textnormal{Pic}(X)$ is injective. The inverse of the matrix $H$ defined above has all positive entries, so each $[\omega_{\mathcal{A}^\vee/X,\tilde\tau}]$ and $[\omega_{\mathcal{A}^\vee/X,\tilde\tau^c}]$  can be written as a linear combination of partial Hasse invariants $[h_{\tilde\tau}]=[\mathcal{O}_X(X_\tau)]$. The strata $X_\tau$'s are disjoint from the cusps, it follows that each $\mathcal{O}_X(X_\tau)$ extends to a line bundle on $X^{\textnormal{min}}$. This shows that all $[\omega_{\tilde\tau}]$ and $[\omega_{\tilde\tau^c}]$ belong to the image of $j_{\mathbb{Q}}$.
\end{proof}

\begin{notation}
We denote by $[\omega_{\mathcal{A}/X,\tau}]$, or simply $[\omega_\tau]$ if no confusion arises, for $\tau\in \Sigma_{\infty,1}$, to be the class $[\omega_{\mathcal{A}/X,\tilde\tau}]=[\omega_{\mathcal{A}/X,\tilde\tau^c}]$ in $\textnormal{Pic}(X)_{\mathbb{Q}}$. According to the lemma, this is well defined and does not depend on the choice of the prime $\mathfrak{q}$ of $\mathcal{O}_E$ above $\mathfrak{p}$. In particular, $[h_\tau]=p^{n_\tau}[\omega_{\sigma^{-n_\tau}\tau}]-[\omega_{\tau}]$ in $\textnormal{Pic}(X)_{\mathbb{Q}}$.
\end{notation}

\begin{lemma}\label{calculation}
We have $p^{n_\tau}[\omega_{\sigma^{-n_\tau}\tau}]+[\omega_{\tau}]=0$ in $\textnormal{Pic}(X_\tau)_{\mathbb{Q}}$.
\begin{proof} Take the lifting $\tilde\tau$ of $\tau$. 
The essential Verschiebung $V_{es,\tilde\tau}^{n_\tau}$ has image $\omega^{(p^{n_\tau})}_{\mathcal{A}^\vee/X,\sigma^{-n_\tau}\tilde\tau}$, so its kernel is a locally free subbundle of $\textnormal{H}^{\textnormal{dR}}_1(\mathcal{A}/X)_{\tilde\tau}$ of rank 1. On $X_\tau$, $V^{n_\tau}_{es,\tilde\tau}$ is zero on $\omega_{\mathcal{A}^\vee/X,\tilde\tau}$, so its kernel is exactly $\omega_{\mathcal{A}^\vee/X,\tilde\tau}$. This could be interpreted as the short exact sequence
\begin{equation}
0\longrightarrow \omega_{\mathcal{A}^\vee/X,\tilde\tau}\longrightarrow \textnormal{H}^{\textnormal{dR}}_1(\mathcal{A}/X)_{\tilde\tau}\xrightarrow{V_{es,\tilde\tau}^{n_\tau}}  \omega^{(p^{n_\tau})}_{\mathcal{A}^\vee/X,\sigma^{-n_\tau}\tilde\tau}\longrightarrow 0.
\end{equation}
So in Pic($X_\tau$), 
\begin{equation}
0=[\wedge^2 \mathrm{H}^{\rm dR}_1(\mathcal{A}/X)_{\tilde\tau}]=[\omega_{\mathcal{A}^\vee/X,\tilde\tau}]+[\omega^{(p^{n_\tau})}_{\mathcal{A}^\vee/X,\sigma^{-n_\tau}\tilde\tau}]=[\omega_\tau]+p^{n_\tau}[\omega_{\sigma^{-n_\tau}\tau}].
\end{equation}
\end{proof}
\end{lemma}

\subsection{Statement of the Main Theorem}

We consider the rational Picard group Pic$(X)_\mathbb{Q}=$Pic$(X)\otimes_\mathbb{Z}\mathbb{Q}$ in the rest of this paper. The notions of ampleness and numerical effectiveness for usual line bundles on varieties extend naturally to the rational Picard group by Chapter 1 in \citep{Lazarsfeld}.

The main theorem of this paper is following.

\begin{theorem}\label{main}
For a tuple of numbers $\underline{t}=(t_\tau)_{\tau\in\Sigma_{\infty,1}}\in \mathbb{Q}^{\Sigma_{\infty,1}}$, the Picard class $[\omega^{\underline{t}}]\stackrel{\text{def}}{=}\sum\limits_{\tau\in\Sigma_{\infty,1}}t_\tau[\omega_\tau]$ is ample if and only if
	\begin{equation}\label{inequality}
		p^{n_\tau}t_\tau>t_{\sigma^{-n_\tau}\tau},\quad \forall\tau\in \Sigma_{\infty,1}.
	\end{equation}
Note that the relation implies that all $t_\tau$'s are necessarily positive.
\end{theorem}

\section{The Goren--Oort Stratification}\label{Stratification}

\subsection{Definition and description of Goren--Oort stratification}

In this section, we introduce a stratification on the special fiber $X$ of the Shimura variety and study its geometry. This will be the key to our proof of the main theorem.

\medskip

\noindent\textbf{Definition}
For any subset $T\subseteq\Sigma_{\infty,1}$, the \emph{Goren--Oort stratum} (or \emph{GO stratum} for short) for $T$ is defined to be the intersection
	\begin{equation}
		X_T\overset{\text{def}}{=}\bigcap_{\tau\in T}X_\tau.
	\end{equation}

The main theorem of this section, roughly speaking, realizes the GO stratum for $T$ as a $\mathbb{P}^1$-power bundle over the special fiber of another unitary Shimura variety. For this, we define a new partition
$\mathcal{P}':\Sigma_\infty=\Sigma'_{\infty,0}\bigsqcup\Sigma'_{\infty,1}\bigsqcup\Sigma'_{\infty,2}$ with
\begin{equation}
\Sigma'_{\infty,0}=\Sigma_{\infty,0}\bigsqcup T'_0,\quad \Sigma'_{\infty,1}=\Sigma_{\infty,1}\backslash T',\quad \Sigma'_{\infty,2}=\Sigma_{\infty,2}\bigsqcup T'_2
\end{equation}
for some $T'=T_0'\sqcup T_2'\subseteq\Sigma_{\infty,1}$. This will give the signature condition for this auxiliary Shimura variety.

\begin{assumption}
Till the end of this paper, we make an additional assumption that $T_{/\mathfrak{p}}\subsetneq \Sigma_{\infty/\mathfrak{p},1}$ for all $\mathfrak{p}\in\Sigma$.
\end{assumption}

Let $T_{/\mathfrak{p}}=T\cap \Sigma_{\infty/\mathfrak{p}}$. We decompose $T_{/\mathfrak{p}}$ into cycles $T_{/\mathfrak{p}}=\bigsqcup_{i} C_i$.
Here by a cycle $C_i$ we mean that $C_i=\{\tau_1,\tau_2=\sigma^{-n_{\tau_1}}\tau_1,\dots,\tau_{m_i}=\sigma^{-n_{\tau_{m_i-1}}}\tau_{m_i-1} \}$ for some $\tau_1\in T_{/\mathfrak{p}}$, such that $\sigma^{-n_{\tau_{m_i}}}\tau_{m_i}\notin T_{/\mathfrak{p}}$ and there does not exist $\tau\in T_{/\mathfrak{p}}$ with $\sigma^{-n_\tau}\tau=\tau_1$. Then we define $C_i'=C_i$ if the cardinality $m_i$ is even, and $C_i'=C_i\bigsqcup \{\sigma^{-n_{\tau_{m_i}}}\tau_{m_i}\}$ if $m_i$ is odd. Finally let $T'_{/\mathfrak{p}}=\bigsqcup_i C'_i$ and $T'=\bigsqcup T'_{/\mathfrak{p}}$. In other words, we only focus on places with signature 1 and extend all cycles to have even cardinality. We also define $I_T=T'\backslash T$ to be the set of additional places we add at the end of each odd cycle.

The construction of $T'_0$ and $T'_2$ is a bit more complicated. For each $C'_i=\{\tau_1,\tau_2=\sigma^{-n_{\tau_1}},\dots,\tau_{2m_i}=\sigma^{n_{\tau_{2m_i-1}}}\tau_{2m_i-1}\}$, decompose it into the disjoint union of $C'_{i,0}=\{\tau_1,\tau_3,\dots,\tau_{2m_i-1}\}$ and $C'_{i,2}=\{\tau_2,\tau_4,\dots,\tau_{2m_i}\}$. Then let $T'_{/\mathfrak{p},\delta}=\bigsqcup_iC'_{i,\delta}$ and $T'_\delta=\bigsqcup T'_{/\mathfrak{p},\delta}$ for $\delta=0$ or $2$. We further introduce subsets 
\begin{equation}
\Delta_i=\bigcup_{j=\text{odd},1\le j\le 2m_i}\{\tau_j,\sigma^{-1}\tau_j,\ldots,\sigma^{-(n_{\tau_j}-1)}\tau_j\}.
\end{equation}
We then define $\Delta(T)_{/\mathfrak{p}}=\cup \Delta_i$ and $\Delta(T)=\cup_{\mathfrak{p}|p}\Delta(T)_{/\mathfrak{p}}$.

According to the previous section, this partition yields a new algebraic group $G_{\mathcal{P}'}$. Note that $\#\Sigma_{\infty,1}-\#\Sigma_{\infty,1}'$ is even. Following \citep{Helm}, we will choose a compact open subgroup $K'\subseteq G_{\mathcal{P}'}(\mathbb{A}^f_\mathbb{Q})$ corresponding to $K\subseteq G_{\mathcal{P}}(\mathbb{A}^f_{\mathbb{Q}})$ later. This will give us the auxiliary Shimura variety $Sh_{K'}(G_{\mathcal{P}'})$.

\begin{example}
Suppose $\Sigma_{\infty/\mathfrak{p}}=\{\tau_1,\tau_2,\dots,\tau_{12}\}$ with the partition $\mathcal{P}$ given by
\begin{equation}
\begin{aligned}
&\Sigma_{\infty,0}=\{\tau_3,\tau_4,\tau_9,\tau_{12}\},\\
&\Sigma_{\infty,1}=\{\tau_1,\tau_2,\tau_5,\tau_7,\tau_{10},\tau_{11}\},\\
&\Sigma_{\infty,2}=\{\tau_6,\tau_8\}.
\end{aligned}
\end{equation}
and 
\begin{equation}
T_{/\mathfrak{p}}=\{\tau_2,\tau_7,\tau_{10}\}.
\end{equation}
We have a decomposition $T_{/\mathfrak{p}}=C_1\sqcup C_2$ with
\begin{equation}
C_1=\{\tau_2\},\quad C_2=\{\tau_7,\tau_{10}\}.
\end{equation}
Correspondingly,
\begin{equation}
C_1'=\{\tau_2,\tau_5\},\quad C_2'=\{\tau_7,\tau_{10}\}.
\end{equation}
Both of them have even cardinality. $T_{/\mathfrak{p}}'=C_1'\sqcup C_2'=\{\tau_2,\tau_5,\tau_7,\tau_{10}\}=T_{/\mathfrak{p},0}'\sqcup T_{/\mathfrak{p},2}'$. Here
\begin{equation}
T_{/\mathfrak{p},0}'=\{\tau_2,\tau_7\},\quad T_{/\mathfrak{p},2}'=\{\tau_5,\tau_{10}\}.
\end{equation}
So in the new decomposition $\mathcal{P}'$, we have
\begin{equation}
\Sigma_{\infty/\mathfrak{p},0}'=\{\tau_2,\tau_3,\tau_4,\tau_7,\tau_9,\tau_{12}\},\quad \Sigma_{\infty/\mathfrak{p},1}'=\{\tau_1,\tau_{11}\},\quad \Sigma_{\infty/\mathfrak{p},2}'=\{\tau_5,\tau_6,\tau_8,\tau_{10}\}.
\end{equation}
with $\Delta(T)=\{\tau_2,\tau_3,\tau_4,\tau_7,\tau_8,\tau_9\}$ and $I_T=\{\tau_5\}$.
\end{example}

\begin{theorem}\label{description}
	The $GO$ stratum $X_T$ is isomorphic to a $(\mathbb{P}^1)^{I_T}$-bundle over $Sh_{K'}(G_{\mathcal{P}'})_k$. 
\end{theorem}

The proof of this theorem is somewhat technical, and was essentially contained in Section 5 of \citep{Tian-Xiao}. We recommend the readers to accept this theorem and jump to Section 3.2 directly when reading this paper the first time.

\begin{proof}
	We will construct a correspondence $X_T\leftarrow Y\rightarrow Z\rightarrow Sh_{K'}(G_{\mathcal{P}'})_k$ and prove that the first two arrows are isomorphisms and the last arrow realizes $Z$ as a $(\mathbb{P}^1)^{I_T}$-bundle over $Sh_{K'}(G_{\mathcal{P}'})_k$.

\medskip
\noindent
\textbf{Step 1} (The level structure) According to \citep{Helm}, Corollary 8.2, since $\#\Sigma_{\infty,1}-\#\Sigma_{\infty,1}'$ is even, there exists a $(V',\langle\cdot,\cdot\rangle')$ isomorphic to $(V,\langle\cdot,\cdot\rangle)$ at all finite places but with signature condition $(r'_\tau,s'_\tau)$ for all infinite place $\tau$ given by the partition $\mathcal{P}'$. This gives an isomorphism between the finite adelic points of algebraic groups $\iota:G_{\mathcal{P}}(\mathbb{A}^f_{\mathbb{Q}})\simeq G_{\mathcal{P}'}(\mathbb{A}^f_{\mathbb{Q}})$. We take $K'=\iota^{-1}(K)\subseteq G_{\mathcal{P}'}(\mathbb{A}^f_\mathbb{Q})$.

\medskip
\noindent
\textbf{Step 2} (Constructing $Y$ and $Z$) We define $Y$ to be the moduli scheme over $k$ which associates to every locally noetherian $k$-scheme $S$ the set of equivalent classes of tuples $(A,\lambda,\eta,B,\lambda',\eta',\phi)$ where
\begin{enumerate}
	\item $(A,\lambda,\eta)\in X_T(S)$;
	\item $(B,\lambda',\eta')\in Sh_{K'}(G_{\mathcal{P}'})(S)$;
	\item $\phi:A\to B$ is an $\mathcal{O}_E$-isogeny such that $p\lambda=\phi^\vee\circ\lambda'\circ\phi$, such that the induced map on de Rham homology
\begin{equation}
\phi_{A,\ast,\tilde\tau}:\text{H}_1^{\text{dR}}(A/S)_{\tilde\tau}\to \text{H}_1^{\text{dR}}(B/S)_{\tilde\tau}
\end{equation}
are isomorphisms for $\tilde\tau\in\Sigma_{E/\infty}$ unless $\tau$ lies in $\Delta(T)$. For $\tau\in\Delta(T)$, we require that
\begin{equation}
\text{Ker}(\phi_{A,\ast,\tilde\tau})=\text{Im}(F^n_{A,es,\tilde\tau}).
\end{equation}

Here, for $\tau\in\Delta(T)$, $n$ is defined to be the unique positive integer such that $\tilde\tau,\sigma^{-1}\tilde\tau,\ldots,\sigma^{1-n}\tilde\tau$ all belong to $\Delta(T)$, but $\sigma^{-n}\tilde\tau$ does not. Note that if $\tau\in T$, then $n=n_\tau$.
	\item The level structures are compatible, i.e., $V^{(p)}(\phi)\circ \eta=\eta'$ as morphisms from $V\otimes_\mathbb{Q}\mathbb{A}^{(p)}_f$ to $V^{(p)}(B)$ modulo $K'$ (identified with $K$ via \textbf{Step 1});
\end{enumerate}
Two tuples $(A,\lambda,\eta,B,\lambda',\eta',\phi)$ and $(\tilde{A},\tilde\lambda,\tilde\eta,\tilde{B},\tilde\lambda',\tilde\eta',\tilde\phi)$ in $Y(S)$ are equivalent if there are $\mathcal{O}_E$-equivariant prime-to-$p$ quasi-isogenies $\psi_A:A\to \tilde{A}$ and $\psi_B:B\to \tilde{B}$ such that
\begin{enumerate}
	\item $\psi_A^\vee\circ\tilde\lambda\circ\psi_A=\lambda$ and $\tilde\eta=\psi_A\circ\eta$, that is, $\psi_A$ gives the equivalence between $(A,\lambda,\eta)$ and $(\tilde{A},\tilde\lambda,\tilde\eta)$ as points of $X(S)$;
	\item $\psi_B^\vee\circ\tilde\lambda'\circ\psi_B=\lambda'$ and $\tilde\eta'=\psi_B\circ\eta$, that is, $\psi_B$ gives the equivalence between $(B,\lambda',\eta')$ and $(\tilde{B},\tilde\lambda',\tilde\eta')$ as points of $Sh_{K'}(G_{\mathcal{P}'})(S)$;
	\item $\tilde\phi\circ\psi_A=\psi_B\circ\phi$.
\end{enumerate}

Then we define the $Z$ to be the moduli scheme over $k$ which associates to any locally noetherian $k$-scheme $S$ the set of isomorphism classes of tuples $(B,\lambda',\eta',\mathcal{J})$ where
\begin{enumerate}
	\item $(B,\lambda',\eta')\in Sh_{K'}(G_{\mathcal{P}'})(S)$;
	\item $\mathcal{J}$ is a collection of line subbundles $J_{\tilde\tau}\subseteq \text{H}_1^{\text{dR}}(B/S)_{\tilde\tau}$ for each $\tau\in I_T$.
\end{enumerate}
Two tuples $(B,\lambda',\eta',\mathcal{J})$ and $(\tilde{B},\tilde\lambda',\tilde\eta',\tilde{\mathcal{J}})$ in $Z(S)$ are equivalent if there is an $\mathcal{O}_E$-equivariant prime-to-$p$ quasi isogeny $\psi_B:B\to \tilde{B}$ such that
\begin{enumerate}
	\item $\psi_B^\vee\circ\tilde\lambda'\circ\psi_B=\lambda'$ and $\tilde\eta'=\psi_B\circ\eta$, that is, $\psi_B$ gives the equivalence between $(B,\lambda',\eta')$ and $(\tilde{B},\tilde\lambda',\tilde\eta')$ as points of $Sh_{K'}(G_{\mathcal{P}'})(S)$;
	\item For each $\tau\in I_T$, $\psi_B^\ast \tilde{J}_{\tilde\tau}=J_{\tilde\tau}$.
\end{enumerate}

It is immediate from the definition that $Z$ is a $(\mathbb{P}^1)^{I_T}$-bundle over $Sh_{K'}(G_{\mathcal{P}'})_k$.

\medskip
\noindent
\textbf{Step 3} (Proof of $X_T\simeq Y$: Bijectivity on closed points) Clearly, there exists a natural forgetful map from $Y\to X$ by sending $(A,\lambda,\eta,B,\lambda',\eta',\phi)\in Y(S)$ to $(A,\lambda,\eta)\in X(S)$. We will use covariant Dieudonn$\acute{\text{e}}$ theory to prove that this map is bijective on geometric closed points and tangent spaces, hence inducing an isomorphism $X_T\simeq Y$.

In the following, we assume $k=\bar{k}$. For $x=(A,\lambda,\eta)\in X_T(k)$, associated to $A[p^\infty]$ is a free $W(k)$-module $\tilde{\mathbb{D}}(A)$ which decomposes into $\tilde{\mathbb{D}}(A)_{\tilde\tau}$ and $\tilde{\mathbb{D}}(A)_{\tilde\tau^c}$-parts via the $\mathcal{O}_E$-action, with a perfect pairing $\tilde{\mathbb{D}}(A)_{\tilde\tau}\times\tilde{\mathbb{D}}(A)_{\tilde\tau^c}\to W(k)$ induced by the polarization $\lambda$. Each summand is a free $W(k)$-module of rank 2. We first construct $B$ from $A$, for this, define $W(k)$-modules $\tilde{\mathbb{D}}(A)_{\tilde\tau}\subseteq M_{\tilde\tau}\subseteq p^{-1}\tilde{\mathbb{D}}(A)_{\tilde\tau}$ by
\begin{equation}
M_{\tilde\tau}=
\begin{cases}
p^{-1}F^n_{A,es}(\tilde{\mathbb{D}}(A)_{\sigma^{-n}\tilde\tau}), & \text{if $\tilde\tau\in\Delta(T)$.}\\
\tilde{\mathbb{D}}(A)_{\tilde\tau},  & \text{if $\tilde\tau\notin\Delta(T)$,}
\end{cases}
\end{equation}
We take $M_{\tilde\tau^c}$ to be the dual of $M_{\tilde\tau}$ under the pairing $\tilde{\mathbb{D}}(A)_{\tilde\tau}\times\tilde{\mathbb{D}}(A)_{\tilde\tau^c}\to W(k)$.

Now we check
\begin{equation}\label{FV stability}
V(M_{\tilde\tau})\subseteq M_{\sigma^{-1}\tilde\tau},\quad F(M_{\sigma^{-1}\tilde\tau})\subseteq M_{\tilde\tau}.
\end{equation}
for all $\tilde\tau$. We distinguish four cases
\begin{itemize}
	\item If $\tau,\sigma^{-1}\tau\notin \Delta(T)$, then $M_{\tilde\tau}=\tilde{\mathbb{D}}(A)_{\tilde\tau}$ and $M_{\sigma^{-1}\tilde\tau}=\tilde{\mathbb{D}}(A)_{\sigma^{-1}\tilde\tau}$, hence (\ref{FV stability}) is clear;
	\item If $\tau\in\Delta(T)$ but $\sigma^{-1}\tau\notin\Delta(T)$, then $M_{\tilde\tau}=p^{-1}F_{A}(\tilde{\mathbb{D}}(A)_{\tilde\tau})$ and $M_{\sigma^{-1}\tilde\tau}=\tilde{\mathbb{D}}(A)_{\sigma^{-1}\tilde\tau}$, hence $F_A(M_{\sigma^{-1}\tilde\tau})\subseteq M_{\tilde\tau}$ and $V_A(M_{\tilde\tau})=M_{\sigma^{-1}\tilde\tau}$;
	\item If $\tau\notin\Delta(T)$ but $\sigma^{-1}\tau\in\Delta(T)$, then it follows from the construction that $\sigma^{-1}\tau\in T$. So $\text{Im}(F^{n_\tau}_{A,es})=\text{Ker}(V^{n_\tau}_{A,es})=\omega_{A^\vee,\tilde\tau}$ by the vanishing of the partial Hasse invariant. In this case, $n=n_\tau$, $M_{\tilde\tau}=\tilde{\mathbb{D}}(A)_{\tilde\tau}$ and $M_{\sigma^{-1}\tilde\tau}=p^{-1}\omega_{A^\vee,\sigma^{-1}\tilde\tau}$. Then $F_A(M_{\sigma^{-1}\tilde\tau})=M_{\tilde\tau}$ and $V(M_{\tilde\tau})=pM_{\sigma^{-1}\tilde\tau}$;
	\item If $\tau,\sigma^{-1}\tau\in\Delta(T)$, suppose $n$ is the integer defined for $\tilde\tau$, then $M_{\tilde\tau}=p^{-1}F_{A,es}^n(\tilde{\mathbb{D}}(A)_{\sigma^{-n}\tilde\tau})$ and $M_{\sigma^{-1}\tilde\tau}=p^{-1}F_{A,es}^{n-1}(\tilde{\mathbb{D}}(A)_{\sigma^{-n}\tilde\tau})$. So $F_A(M_{\sigma^{-1}\tilde\tau})=M_{\tilde\tau}$ and $V_A(M_{\tilde\tau})=pM_{\sigma^{-1}\tilde\tau}$.
\end{itemize}

Consequently, if we define $M=\bigoplus_{\tau\in\Sigma_{\infty}}(M_{\tilde\tau}\oplus M_{\tilde\tau^c})$, then $M$ is a Dieudonn$\acute{\text{e}}$ module such that $\tilde{\mathbb{D}}(A)\subseteq M\subseteq p^{-1}\tilde{\mathbb{D}}(A)$, with induced maps $F$ and $V$ on $M$. Consider the quotient $M/\tilde{\mathbb{D}}(A)$, according to Dieudonn$\acute{\text{e}}$ theory, it corresponds to a finite subgroup scheme $K\subseteq A[p]$ stable under the action of $\mathcal{O}_E$. We then put $B=A/K$, equipped with a map $\iota_B:\mathcal{O}_E\to \text{End}(B)$, and with a natural quotient morphism $\phi:A\to B$. The induced morphism on Dieudonn$\acute{\text{e}}$ modules corresponds to $\tilde{\mathbb{D}}(A)\subseteq M$.

We define the quasi-polarization of $B$ via
\begin{equation}
\lambda_B:B\xleftarrow{\phi} A\stackrel{\lambda}{\dashrightarrow} A^\vee\xleftarrow{\phi^\vee}B^\vee.
\end{equation}

It remains to check that $B$ satisfies the desired signature condition. The point is that we always have the diagram
\begin{equation}
\begin{tikzcd}
{\mathbb{D}}(A)_{\tilde\tau} \arrow[r, "V"] \arrow[d, "{\phi_{\ast,\tilde\tau}}"'] & {\mathbb{D}}(A)_{\sigma^{-1}\tilde\tau} \arrow[d, "{\phi_{\ast,\sigma^{-1}\tilde\tau}}"]  \\
{\mathbb{D}}(B)_{\tilde\tau} \arrow[r, "V"]  & {\mathbb{D}}(B)_{\sigma^{-1}\tilde\tau}
\end{tikzcd}.
\end{equation}
with $\mathbb{D}(A)_{\tilde\tau}=\tilde{\mathbb{D}}(A)_{\tilde\tau}/p\tilde{\mathbb{D}}(A)_{\tilde\tau}$ and similarly for $B$. Then the arrows are maps of $2$-dimensional $k$-vector spaces. By construction, we have the dimensions of the kernels and images of the vertical arrows; we also know that dim(Im $V_A$)=$2-r_{\sigma^{-1}\tau}=s_{\sigma^{-1}\tau}$. We again distinguish four cases as before:
\begin{itemize}
	\item If $\tau,\sigma^{-1}\tau\notin \Delta(T)$, then the vertical arrows are isomorphisms, so $s'_{\sigma^{-1}\tau}=s_{\sigma^{-1}\tau}$;
	\item If $\tau\in\Delta(T)$ but $\sigma^{-1}\tau\notin\Delta(T)$, then $\phi_{\ast,\tilde\tau}$ has 1-dimensional kernel and $\phi_{\ast,\tilde\tau}$ is an isomorphism, and $s'_{\sigma^{-1}\tau}=s_{\sigma^{-1}\tau}-1=0$ since $\sigma^{-1}\tilde\tau\in \Sigma_{\infty,1}$;
	\item If $\tau\notin\Delta(T)$ but $\sigma^{-1}\tau\in\Delta(T)$, then $\phi_{\ast,\tilde\tau}$ is an isomorphism and $\phi_{\ast,\sigma^{-1}\tilde\tau}$ has 1 dimensional kernel, and $s'_{\sigma^{-1}\tau}=s_{\sigma^{-1}\tau}+1=2$ since $\sigma^{-1}\tilde\tau\in \Sigma_{\infty,1}$;
	\item If $\tau,\sigma^{-1}\tau\in\Delta(T)$, then both $\phi_{\ast,\tilde\tau}$ and $\phi_{\ast,\sigma^{-1}\tilde\tau}$ have 1-dimensional kernels, so $s'_{\sigma^{-1}\tau}=s_{\sigma^{-1}\tau}$.
\end{itemize}
This corresponds to the signature condition in the moduli problem for $B$.

The quasi-polarization induces isomorphisms of prime-to-$p$ Tate modules. Via the choice of $K'$ in \textbf{Step 1}, the $K$-level structure of $A$ corresponds to a $K'$-level structure of $B$.

Thus we have construct a bijection on closed points of the moduli schemes $X_T$ and $Y$.

\medskip
\noindent\textbf{Step 4} (Proof of $X_T\simeq Y$: Bijectivity on tangent spaces) Suppose that $x=(A,\lambda,\eta)$ corresponds to $y=(A,\lambda,\eta,B,\lambda_B,\eta_B,\phi)$ under the bijection in \textbf{Step 3}. We recall the following theory of deformation of abelian schemes:

Suppose $S$ is a nilpotent thickening of $S_0$ equipped with divided power structure. A basic example used in the following is $\text{Spec}\ k\hookrightarrow \text{Spec}\ k[\varepsilon]/(\varepsilon^2)$ with $k$ a perfect field of characteristic $p$. For an abelian scheme $A_0/S_0$, let $\text{H}^{\text{cris}}_1(A_0/S_0)_S=\text{H}^1_{\text{cris}}(A^\vee_0/S_0)_S$ be the relative crystalline homology group. There is a canonical isomorphism
\begin{equation}
\text{H}_1^{\text{cris}}(A_0/S_0)_S\otimes_{\mathcal{O}_S}\mathcal{O}_{S_0}\simeq \text{H}_1^{\text{dR}}(A_0/S_0).
\end{equation}
If we have moreover an abelian scheme $A/S$ with $A\times_S S_0\simeq A_0$, then there is a canonical Hodge filtration
\begin{equation}
0\longrightarrow \omega_{A^\vee/S}\longrightarrow \text{H}_1^{\text{cris}}(A_0/S_0)_S\longrightarrow \text{Lie}(A/S)\longrightarrow 0.
\end{equation}
Hence, $\omega_{A^\vee/S}$ is a local direct factor of $\text{H}^{\text{cris}}_1(A_0/S_0)_S$ that lifts the subbundle $\omega_{A_0^\vee/S_0}$. If we denote by $AV_S$ to be the category of abelian schemes over $S$ and $AV_{S_0}^+$ the category of pairs $(A_0,\omega_S)$, where $A_0\in AV_{S_0}$ and $\omega_S\subseteq \text{H}_1^{\text{cris}}(A_0/S_0)_S$ is a line subbundle with $\omega_S\otimes_{\mathcal{O}_{S}}\mathcal{O}_{S_0}=\omega_{A_0^\vee/S_0}$. We thus obtain a functor
\begin{equation}
\begin{aligned}
AV_{S}&\longrightarrow AV_{S_0}^+\\
A&\longrightarrow (A\times_S S_0,\omega_{A^\vee/S}).
\end{aligned}
\end{equation}
\begin{theorem}\label{deformation}
The above functor is an equivalence of categories. Moreover,
\begin{enumerate}
	\item If $A_0$ is equipped with an action of a (not necessarily commutative) ring $\iota:R\to \textnormal{End}(A_0)$, then the action lifts to $A$ if and only if $\omega_{A^\vee/S}\subseteq\textnormal{H}_1^{\textnormal{dR}}(A_0/S_0)_S$ is $R$-stable.
	\item If $A_0$ is equipped with a quasi-polarization $\lambda_0:A_0\to A_0^\vee$, then $\lambda_0$ induces an alternating pairing $\langle\cdot,\cdot\rangle_{\lambda_0}$ on $\textnormal{H}_1^{\textnormal{cris}}(A_0/S_0)_S$ by \textnormal{\citep{BBM}}, which is perfect if $\lambda_0$ is prime-to-$p$. Then there exists a (necessarily unique) quasi-polarization $\lambda:A\to A^\vee$ that lifts $\lambda_0$ if and only if $\omega_S$ is isotropic under $\langle\cdot,\cdot\rangle_{\lambda_0}$.
\end{enumerate}
\end{theorem}
In other words, lifting an abelian scheme (with extra structures) over $S_0$ to $S$ is equivalent to lifting the relative differential sheaf to a (``stable'') subbundle of the crystalline homology group.

Back to the morphism $Y\to X_T$ with $y=(A,\lambda,\eta,B,\lambda',\eta',\phi)$ mapped to $x=(A,\lambda,\eta)$. Write $\mathbb{I}=\text{Spec}\  k[\varepsilon]/(\varepsilon^2)$. A tangent vector of $X_T$ at $x$ is the tuple $x_{\mathbb{I}}=(A_{\mathbb{I}},\lambda_{A,\mathbb{I}},\eta_{A,\mathbb{I}})$. Giving the abelian scheme $A_{\mathbb{I}}$ with $\mathcal{O}_E$-action lifting $(A,\lambda)$ is equivalent to lifting all $\omega_{A^\vee/k,\tilde\tau}$ and $\omega_{A^\vee/k,\tilde\tau^c}$ to subbundles $\omega_{A^\vee,\mathbb{I},\tilde\tau},\omega_{A^\vee,\mathbb{I},\tilde\tau^c}$ of $\text{H}^{\text{cris}}_1(A/k)_{\mathbb{I},\tilde\tau}$ and $\text{H}^{\text{cris}}_1(A/k)_{\mathbb{I},\tilde\tau^c}$ respectively; and to lift the quasi-polarization, we need to require that $\omega_{A^\vee,\mathbb{I},\tilde\tau}$ and $\omega_{A^\vee,\mathbb{I},\tilde\tau^c}$ are annilators of each other. Note that if $\tau\in \Sigma_{\infty,0}$ or $\tau\in\Sigma_{\infty,2}$, the subbundles are either $0$ or the whole bundle, so there is a unique lifting for these places. It remains to lift the line subbundles $\omega_{A^\vee/S,\tilde\tau}$. Moreover, the condition that partial Hasse invariants $h_{\tilde\tau}$ for $\tau\in T$ vanishes is equivalent to that the liftings $\omega_{A^\vee,\mathbb{I},\tilde\tau}=F^{n_\tau}_{A,es}(\text{H}_1^{\text{cris}}(A^{(p^{n_\tau})}/k)_{\mathbb{I},\tilde\tau})$ for $\tau\in T$, where the Frobenius and Verschiebung, essential Frobenius and essential Verschiebung maps on crystalline homology groups are defined in the same manner as on de Rham homology.

Starting from $x_{\mathbb{I}}$, we will show that there exists a unique deformation $(B_{\mathbb{I}},\lambda_{B,\mathbb{I}},\eta_{B,\mathbb{I}},\phi_{\mathbb{I}})$ of $(B,\lambda_B,\eta_B,\phi)$ such that $(A_{\mathbb{I}},\lambda_{A,\mathbb{I}},\eta_{A,\mathbb{I}},B_{\mathbb{I}},\lambda_{B,\mathbb{I}},\eta_{B,\mathbb{I}},\phi_{\mathbb{I}})\in Y(\mathbb{I})$.

We first construct $B_{\mathbb{I}}$, which is a lifting of $B$.
\begin{itemize}
	\item If either $\tau$ or $\sigma\tau$ belongs to ${\Delta}(T)$, the signature condition shows that $\omega_{B^\vee,\tilde\tau}$ is either 0 or the whole $\text{H}_1^{\text{dR}}(B/k)_{\tilde\tau}$, so there is a unique lifting;
	\item If $\tau,\sigma\tau\notin\Delta(T)$, we have isomorphisms in the construction of $B$:
\begin{equation}
\begin{aligned}
&\phi_{A,\ast,\tilde\tau}:\text{H}_1^{\text{dR}}(A/k)_{\tilde\tau}\longrightarrow \text{H}^{\text{dR}}_1(B/k)_{\tilde\tau},\\
&\phi_{A,\ast,\sigma\tilde\tau}:\text{H}_1^{\text{dR}}(A/k)_{\sigma\tilde\tau}\longrightarrow \text{H}^{\text{dR}}_1(B/k)_{\sigma\tilde\tau}.
\end{aligned}
\end{equation}
We take $\omega_{B^\vee,\mathbb{I},\tilde\tau}\subseteq \text{H}_1^{\text{cris}}(B/k)_{\mathbb{I},\tilde\tau}$ to be the image of $\omega_{A^\vee,\mathbb{I},\tilde\tau}\subseteq \text{H}_1^{\text{cris}}(A/k)_{\mathbb{I},\tilde\tau}$ under the morphism $\phi^{\text{cris}}_{\ast,\tilde\tau}$ on crystalline homology. We take the $\tilde\tau^c$-part to be the dual. Then $B_{\mathbb{I}}$ is an abelian scheme equipped with an $\mathcal{O}_E$-action by the construction, and $\phi_{A,\mathbb{I}}:A_{\mathbb{I}}\to B_{\mathbb{I}}$ is $\mathcal{O}_E$-equivariant by \citep{Lan}.
\end{itemize}

We have to check the conditions on $\phi_{\mathbb{I},\ast,\tilde\tau}:\text{H}_1^{\text{dR}}(A_{\mathbb{I}}/\mathbb{I})_{\tilde\tau}\to\text{H}_1^{\text{dR}}(B_{\mathbb{I}}/\mathbb{I})_{\tilde\tau}$, which is canonically identified with
\begin{equation}
\phi_{\mathbb{I},\ast,\tilde\tau}^{\text{cris}}:\text{H}_1^{\text{cris}}(A/k)_{\mathbb{I},\tilde\tau}\to\text{H}_1^{\text{cris}}(B/k)_{\mathbb{I},\tilde\tau}.
\end{equation}
For $\tau\notin \Delta(T)$, the map is an isomorphism since it is for $\phi_{\ast,\tilde\tau}$; for $\tau\in\Delta(T)$, the Frobenius map on $k[\varepsilon]/(\varepsilon^2)$ factors as
\begin{equation}
k[\varepsilon]/(\varepsilon^2)\twoheadrightarrow k\xrightarrow{x\mapsto x^p}k\hookrightarrow k[\varepsilon]/(\varepsilon^2).
\end{equation}
So we have a canonical isomorphism
\begin{equation}
F^n_{A,es,\tilde\tau}(\text{H}_1^{\text{dR}}(A^{(p^n)}_{\mathbb{I}}/\mathbb{I})_{\tilde\tau})=F^n_{A,es,\tilde\tau}(\text{H}_1^{\text{dR}}(A^{(p^n)}/k)_{\tilde\tau})\otimes_kk[\varepsilon]/(\varepsilon^2).
\end{equation}
This shows that $\text{Ker}(\phi_{\mathbb{I},\ast,\tilde\tau}:\text{H}_1^{\text{dR}}(A_{\mathbb{I}}/\mathbb{I})_{\tilde\tau}\to \text{H}_1^{\text{dR}}(B_{\mathbb{I}}/\mathbb{I})_{\tilde\tau})=F_{A,es,\tilde\tau}^n(\text{H}_1^{\text{dR}}(A_{\mathbb{I}}^{(p^n)}/\mathbb{I})_{\tilde\tau})$. Again we take the lift $\omega_{A^\vee,\mathbb{I},\tilde\tau^c}$ to be the annihilator of $\omega_{A^\vee,\mathbb{I},\tilde\tau}$.
Therefore $\phi_{A,\mathbb{I}}$ satisfies the condition in the moduli problem, and the polarization $\lambda_B$ deforms to $\lambda_{\mathbb{I}}$ by construction.

The prime-to-$p$ isogeny $\phi_{\mathbb{I}}:A_{\mathbb{I}}\to B_{\mathbb{I}}$ induces an isomorphism between prime-to-$p$ Tate modules. Via the identification in \textbf{Step 1}, we take the level structure $\eta_{A,\mathbb{I}}$ to be the one corresponding to $\eta_{A,\mathbb{I}}$.

We have shown the bijection on tangent spaces. So we get the isomorphism $X_T\simeq Y$.

\medskip
\noindent
\textbf{Step 5} (Proof of $Y\simeq Z$: Bijectivity on closed points)

We first define a morphism $\pi_2:Y\to Z$. For a locally noetherian $k$-scheme $S$ and $y=(A,\lambda_A,\eta_A,B,\lambda_B,\eta_B,\phi)\in Y(S)$, let $\pi_2(x)=(B,\lambda_B,\eta_B,\mathcal{J})$ where
\begin{equation}
J_{\tilde\tau}=\phi_{\ast,\tilde\tau}(\omega_{A^\vee,\tilde\tau}).
\end{equation}
for $\tau\in I_T$. Note that $\phi_{\ast,\tilde\tau}$ is an isomorphism by construction.

Now suppose $k=\bar{k}$ and $z=(B,\lambda_B,\eta_B,\mathcal{J})\in Z(k)$ is a closed point, we first construct $A$ from $z$. For this, we will construct modules $\tilde{\mathbb{D}}(A)_{\tilde\tau}\subseteq M_{\tilde\tau}\subseteq p^{-1}\tilde{\mathbb{D}}(A)_{\tilde\tau}$. For $\tau\in I_T$, let $\tilde{J}_{\tilde\tau}$ be the inverse image of $\bar{J}_{\tilde\tau}\subseteq \mathbb{D}(B)_{\tilde\tau}$ via the modulo $p$ map $\tilde{\mathbb{D}}(B)\to\mathbb{D}(B)$.
\begin{itemize}
	\item If $\tau\notin \Delta(T)$, then we directly take $M_{\tilde\tau}=p^{-1}\tilde{\mathbb{D}}(B)_{\tilde\tau}$;
	\item If $\tau\in\Delta(T)$, recall the definition of $\Delta(T)$, we assume $\tau=\sigma^{-l}\tau_i$ for some signature 1 place $\tau_i$ in an odd chain $C=\{\tau_1,\tau_2=\sigma^{-n_{\tau_1}}\tau_1,\ldots,\tau_{2m-1}\}$ with $0\le l< n_{\tau_i},1<i<2m-1$. Let $\tau_{2m}=\sigma^{-n_{\tau_{2m-1}}}\tau_{2m-1}\in I_T$ and suppose $\tau=\sigma^a\tau_{2m}$. The composition
\begin{equation}
F^a_{B,es}:\tilde{\mathbb{D}}(B)_{\tilde\tau_{2m}}\xrightarrow[\simeq]{F_{B,es}}\tilde{\mathbb{D}}(B)_{\sigma\tilde\tau_{2m}}\xrightarrow[\simeq]{F_{B,es}}\cdots\xrightarrow[\simeq]{F_{B_es}}\tilde{\mathbb{D}}(B)_{\tilde\tau}.
\end{equation}
is an isomorphism since in this case, $s_{\tau_1},s_{\sigma^{-1}\tau_1},\ldots,s_{\tau_{2m}}=0$ or $2$, so each essential Frobenius $F_{B,es}$ is an isomorphism. We set $M_{\tilde\tau}=p^{-1}F_{B,es}^a(\tilde{J}_{\tilde\tau_{2m}})$.
	\item If $\tau\in\Delta(T)$ but $\tau=\sigma^{-l}\tau_i$ for some signature 1 place $\tau_i$ in an even chain $C=\{\tau_1,\tau_2=\sigma^{-n_{\tau_1}}\tau_1,\ldots,\tau_{2m}\}$ with $0\le l<n_{\tau_i},1<i<2m$. There is no $J$ involved in the construction of this part. Again write $\tau_{2m+1}=\sigma^{n_{\tau_{2m}}}\tau_{2m}$ and $\tau=\sigma^b\tau_{2m+1}$, the morphisms in the composition
\begin{equation}
F^b_{B,es}:\tilde{\mathbb{D}}(B)_{\tilde\tau_{2m+1}}\xrightarrow{F_{B}}\tilde{\mathbb{D}}(B)_{\sigma\tilde\tau_{2m+1}}\xrightarrow[\simeq]{F_{B,es}}\cdots\xrightarrow[\simeq]{F_{B_es}}\tilde{\mathbb{D}}(B)_{\tilde\tau}.
\end{equation}
are isomorphisms except the first one. We put $M_{\tilde\tau}=p^{-1}F^b_{B,es}(\tilde{\mathbb{D}}(B)_{\tilde\tau_{2m+1}})$. The signature condition $s'_{\tau_{2m+1}}=1$ shows that 
\begin{equation}
\text{dim}_k\ \text{Ker}(F_B:\tilde{\mathbb{D}}(B)_{\tilde\tau_{2m+1}}\to \tilde{\mathbb{D}}(B)_{\sigma\tilde\tau_{2m+1}})=1.
\end{equation}
So $\text{dim}_k(M_{\tilde\tau}/\tilde{\mathbb{D}}(B)_{\tilde\tau})=1$.
\end{itemize}

One can check that $M=\bigoplus_{\tau\in\Sigma_{\infty}}(M_{\tilde\tau}\oplus M_{\tilde\tau^c})$ is stable under $F$ and $V$:
\begin{itemize}
	\item If $\tau,\sigma\tau\notin\Delta(T)$, then $M_{\tilde\tau}=p^{-1}\tilde{\mathbb{D}}(B)_{\tilde\tau}$ and $M_{\sigma\tilde\tau}=p^{-1}\tilde{\mathbb{D}}(B)_{\sigma\tilde\tau}$, $F,V$-stability is direct;
	\item If $\sigma\tau\notin\Delta(T)$ but $\tau\in\Delta(T)$, then $M_{\sigma\tilde\tau}=p^{-1}\tilde{\mathbb{D}}(B)_{\sigma\tilde\tau}$ and $M_{\tilde\tau}$ is either $p^{-1}F^a_{B,es}(\tilde{J}_{\tau_{2m}})$ or $p^{-1}F^b_{B,es}(\tilde{\mathbb{D}})_{\tilde\tau_{2m+1}}$. According to the signature condition of $B$, $s'_{\tau}=2$, so $V\tilde{\mathbb{D}}(B)_{\sigma\tilde\tau}=p^{-1}\tilde{\mathbb{D}}(B)_{\tilde\tau}$. Then $VM_{\sigma\tilde\tau}=\tilde{\mathbb{D}}(B)_{\tilde\tau}\subseteq M_{\tilde\tau}$ and $FM_{\tilde\tau}\subseteq p^{-1}\tilde{\mathbb{D}}(B)_{\sigma\tau}=M_{\sigma\tilde\tau}$;
	\item If $\sigma\tau\in\Delta(T)$ but $\tau\notin\Delta(T)$, then either $\tau\in I_T$ or $\tau=\tau_{2m}$ in the above definition. Either $M_{\sigma\tilde\tau}=p^{-1}F_{B,es}(\tilde{J}_{\tau_{2m}})$ or $M_{\sigma\tilde\tau}=p^{-1}F_{B,es}^b(\tilde{\mathbb{D}}(B)_{\tilde\tau_{2m+1}})$, $M_{\tilde\tau}=p^{-1}\tilde{\mathbb{D}}(B)_{\tilde\tau}$. $F,V$-stability follows from the signature condition of $B$ that $s'_{\sigma\tau}=2$.
	\item If $\tau,\sigma\tau\in\Delta(T)$, either $M_{\tilde\tau}=p^{-1}F^a_{B,es}(\tilde{J}_{\tilde\tau_{2m}})$ and $M_{\sigma\tilde\tau}=p^{-1}F^{a+1}_{B,es}(\tilde{J}_{\tilde\tau_{2m}})$ or $M_{\tilde\tau}=p^{-1}F^b_{B,es}(\tilde{\mathbb{D}}(B)_{\tilde\tau_{2m+1}})$ and $M_{\sigma\tilde\tau}=p^{-1}F^{b+1}_{B,es}(\tilde{\mathbb{D}}(B)_{\tilde\tau_{2m+1}})$. In the former case, $F,V$ stability follows from $s'_{\tau}$ and $s'_{\sigma\tau}$ are either 0 or 2. The latter case is direct.
\end{itemize}

We take $M_{\tilde\tau^c}$ to be the dual of $M_{\tilde\tau}$, viewed as lattices in $\tilde{\mathbb{D}}(B)[1/p]_{\tilde\tau}$ and $\tilde{\mathbb{D}}(B)[1/p]_{\tilde\tau^c}$ respectively. It follows that if we take 
\begin{equation}
M=\bigoplus_{\tau\in\Sigma_{\infty}}(M_{\tilde\tau}\oplus M_{\tilde\tau^c}).
\end{equation}
Then $M$ is a Dieudonn$\acute{\text{e}}$ module equipped with an $\mathcal{O}_E$-action. The quotient
\begin{equation}
M/\tilde{\mathbb{D}}(B)\subseteq p^{-1}\tilde{\mathbb{D}}(B)/\tilde{\mathbb{D}}(B).
\end{equation}
corresponds to a finite group scheme $G\subseteq B[p]$ stable under $\mathcal{O}_E$-action. Put $A=B/G$ and let $\psi:B\to A$ be the natural projection. We take $\phi:A\to B$ such that $\phi\circ\psi=[p]_B$.

We define the prime-to-$p$ quasi-polarization of $A$ to be $1/p$ times the composition
\begin{equation}
A\xrightarrow{\phi}B\stackrel{\lambda_B}{\dashrightarrow}B^\vee \xrightarrow{\phi^\vee} A^\vee.
\end{equation}
Here, the middle quasi-isogeny is prime-to-$p$, and $\text{Ker}\ \phi$ is a maximal isotropic subgroup of $A[p]$, so the composition can be written as $\lambda_A\circ p$ with $\lambda_A$ prime-to-$p$.

Via the identification in \textbf{Step 1}, we choose the level structure $\eta_A$ of $A$ corresponding to $\eta_B$.

By a case by case discussion similar as before, we can show that $A$ satisfies the desired signature condition. So we obtain a unique point $y=(A,\lambda_A,\eta_A,B,\lambda_B,\eta_B,\phi)\in Y(k)$ such that $\pi_2(y)=z$. This proves the bijiectivity on closed points.

\medskip
\noindent
\textbf{Step 6} (Proof of $Y\simeq Z$: Bijectivity on tangent spaces)

Continue to suppose that $y=(A,\lambda_A,\eta_A,B,\lambda_B,\eta_B,\phi)\in Y(k)$ is mapped to $z=(B,\lambda_B,\eta_B,\mathcal{J})$ via $\pi_2$. Let $\mathbb{I}=k[\varepsilon]/(\varepsilon^2)$. Suppose $z_{\mathbb{I}}=(B_{\mathbb{I}},\lambda_{B,\mathbb{I}},\eta_{B,\mathbb{I}},\mathcal{J}_{\mathbb{I}})\in Z(\mathbb{I})$, where $\mathcal{J}_{\mathbb{I}}$ is a collection of line subbundles $J_{\mathbb{I},\tilde\tau}\subseteq \text{H}_1^{\text{dR}}(B_{\mathbb{I}}/\mathbb{I})_{\tilde\tau}=\text{H}_1^{\text{cris}}(B/k)_{\mathbb{I},\tilde\tau}$ lifting $J_{\tilde\tau}$ for each $\tau\in I_T$. We have to show that $z_{\mathbb{I}}$ corresponds uniquely to a deformation $y_{\mathbb{I}}$ of $y$ such that $\pi_{\mathbb{I},2}(y_{\mathbb{I}})=z_{\mathbb{I}}$.

We construct the liftings for $\omega_{A^\vee/k,\tilde\tau}$ for all $\tau\in\Sigma_{\infty}$ as follows:
\begin{itemize}
	\item If $\tau,\sigma\tau\in\Delta(T)$, then we have isomorphisms
\begin{equation}
\begin{aligned}
&\phi_{\mathbb{I},\ast,\tilde\tau}^{\text{cris}}:\text{H}_1^{\text{cris}}(A/k)_{\mathbb{I},\tilde\tau}\xrightarrow{\simeq} \text{H}^{\text{cris}}_1(B/k)_{\mathbb{I},\tilde\tau},\\
&\phi_{\mathbb{I},\ast,\sigma\tilde\tau}^{\text{cris}}:\text{H}_1^{\text{cris}}(A/k)_{\mathbb{I},\sigma\tilde\tau}\xrightarrow{\simeq} \text{H}^{\text{cris}}_1(B/k)_{\mathbb{I},\sigma\tilde\tau}.
\end{aligned}
\end{equation}
We then take $\omega_{A^\vee,\mathbb{I},\tilde\tau}=(\phi^{\text{cris}}_{\mathbb{I},\ast,\tilde\tau})^{-1}(\omega_{B^\vee_{\mathbb{I}}/\mathbb{I},\tilde\tau})$.
	\item If $\tau\in I_T$, $\phi_{\mathbb{I},\ast,\tilde\tau}^{\text{cris}}$ is an isomorphism. We take $\omega_{A^\vee,\mathbb{I},\tilde\tau}=(\phi^{\text{cris}}_{\mathbb{I},\ast,\tilde\tau})^{-1}(J_{\mathbb{I},\tilde\tau})$.
	\item If $\tau,\sigma\tau\in \Delta(T)$, then the signature $s_\tau=0$ or $2$, so $\omega_{A^\vee,\tilde\tau}$ is either empty or $\text{H}_1^{\text{dR}}(A/k)$, which admits a unique lifting to $\omega_{A^\vee,{\mathbb{I}},\tilde\tau}$.
	\item For the rest cases, $\tau\in T$, the partial Hasse invariant $h_{\tilde\tau}$ vanishes, so the lifting is completely determined by
\begin{equation}
\omega_{A^\vee,\mathbb{I},\tilde\tau}=\text{Ker}(V^{n_\tau}_{A_\mathbb{I},es,\tilde\tau})=F^{n_\tau}_{A_{\mathbb{I}},es,\tilde\tau}(\text{H}_1^{\text{cris}}(A^{(p^{n_\tau})}/k)_{\mathbb{I},\tilde\tau}).
\end{equation}
\end{itemize}
We take $\omega_{A^\vee,\mathbb{I},\tilde\tau^c}$ to be the orthogonal complement of $\omega_{A^\vee,\mathbb{I},\tilde\tau}$ under the perfect pairing
\begin{equation}
\langle\cdot,\cdot\rangle_{\lambda_A}:\text{H}_1^{\text{cris}}(A/k)_{\mathbb{I},\tilde\tau}\times \text{H}_1^{\text{cris}}(A/k)_{\mathbb{I},\tilde\tau^c}\to k[\varepsilon]/(\varepsilon^2).
\end{equation}
Now according to Theorem $\ref{deformation}$, we obtain an abelian scheme $A_{\mathbb{I}}/\mathbb{I}$ with $\mathcal{O}_E$-action and prime-to-$p$ quasi-polarization $\lambda_{A,\mathbb{I}}$ which lifts $(A,\lambda_A)$. Again via the identification in \textbf{Step 1}, we take the level structure $\eta_{A,\mathbb{I}}$ of $A_{\mathbb{I}}$ corresponding to that of $B_{\mathbb{I}}$. So we get $y_{\mathbb{I}}=(A_{\mathbb{I}},\lambda_{A,\mathbb{I}},\eta_{A,\mathbb{I}},B_{\mathbb{I}},\eta_{B,\mathbb{I}},\lambda_{B,\mathbb{I}},\eta_{B,\mathbb{I}},\phi_{\mathbb{I}})\in Y(\mathbb{I})$ lifting $y$ and $\pi_{\mathbb{I},2}(y_{\mathbb{I}})=z_{\mathbb{I}}$.

Combining the above steps, we have proved the description of Goren-Oort stratum.

\end{proof}

\begin{remark}\label{identifications}
According to the proof, if we let $j:X_T\hookrightarrow X$ be the natural embedding, $\mathcal{A}'$ be the universal abelian scheme over $Sh_{K'}(G_{\mathcal{P}'})_k$ and $\pi:X_T\rightarrow Sh_{K'}(G_{\mathcal{P}'})_k$ be the bundle morphism, then we have identifications $j^\ast\omega_{\mathcal{A}^\vee/X,\tilde\tau}\simeq \pi^\ast\omega_{\mathcal{A}'^\vee/Sh_{K'}(G_{\mathcal{P}'})_k,\tilde\tau}$ and $j^\ast\omega_{\mathcal{A}^\vee/X,\tilde\tau^c}\simeq \pi^\ast\omega_{\mathcal{A}'^\vee/Sh_{K'}(G_{\mathcal{P}'})_k,\tilde\tau^c}$ for all $\tau\in\Sigma'_{\infty,1}$.
\end{remark}

\subsection{Positivity of line bundles}

We collect some tools for the study of ampleness of line bundles in this section. The main reference for this subsection is Section 1.4 of \citep{Lazarsfeld}. Let $X$ be a complete variety over an algebraically closed field $k$ and $\mathcal{L}$ a line bundle on $X$. 

\medskip
\noindent\textbf{Definition}
	A Weil divisor $D$ is called \emph{numerically effective}, or \emph{nef} for short, if the intersection number $(D\cdot C)$ is nonnegative for all irreducible curves $C\subseteq X$; Two divisors $D_1,\ D_2$ are said to be \emph{numerically equivalent}, written $D_1\equiv_{num} D_2$, if $(D_1\cdot C)=(D_2\cdot C)$ for all irreducible curves $C\subseteq X$. 

\medskip
\noindent\textbf{Remark}
It is shown that if two divisors are linear equivalent, then they are numerically equivalent. So nefness can be defined for the elements in the rational Picard group. A line bundle is called \emph{nef} if its corresponding divisor class is nef.

The key to prove Theorem \ref{main} is to prove the following nef version.
\begin{theorem}\label{variant}
	For a tuple of numbers $\underline{t}=(t_\tau)_{\tau\in\Sigma_{\infty,1}}\in \mathbb{Q}^{\Sigma_{\infty,1}}$, the Picard class $[\omega^{\underline{t}}]\stackrel{\text{def}}{=}\sum\limits_{\tau\in\Sigma_{\infty,1}}t_\tau[\omega_\tau]$ is nef if
	\begin{equation}\label{inequality'}
		p^{n_\tau}t_\tau\ge t_{\sigma^{-n_\tau}\tau},\quad \forall\tau\in \Sigma_{\infty,1}.
	\end{equation}
\end{theorem}

\begin{proposition}\label{ampleness}
Let $D,E$ be rational divisors on $X$. If $D$ is ample and $E$ is nef, then their sum $D+E$ is ample. \emph{(\citep{Lazarsfeld}, Corollary 1.4.10)}
\end{proposition}

\begin{proposition}\label{pullback}
Let $X$ be a complete variety and $\mathcal{L}$ a line bundle over $X$. Suppose $f:Y\to X$ is a proper morphism. If $\mathcal{L}$ is nef, then $f^\ast(\mathcal{L})$ is nef. In particular, restrictions of nef bundles to subschemes remains nef. \emph{(\citep{Lazarsfeld}, Example 1.4.4)}
\end{proposition}

\subsection{Necessity of the ampleness criterion}

In this subsection, we prove the necessity part of the main theorem. This argument is essentially due to \citep[\S6]{Tian-Xiao}. So now suppose $X=Sh_K(G_\mathcal{P})_k$ is the special fiber of our given unitary Shimura variety and $[\omega^{\underline{t}}]=\sum_{\tau\in\Sigma_{\infty,1}}t_\tau[\omega_\tau]$ is an ample class over $X$
associated to a tuple of numbers $\underline{t}=(t_\tau)_{\tau\in \Sigma_{\infty,1}}$. We have to show the inequality $p^{n_\tau}t_\tau>t_{\sigma^{-n_\tau}\tau}$ for each $\tau$.

Suppose $\mathfrak{p}$ is a prime above $p$. We distinguish two cases:

\textbf{Case 1}: $\Sigma_{\infty/\mathfrak{p},1}=\{\tau\}$ for some $\tau$. We need to show that $t_{\tau}>0$. Take $T=\Sigma_{\infty,1}-\{\tau\}$ and consider the corresponding GO stratum $X_T$. This is a curve on $X$. For any $\mathfrak{q}\ne\mathfrak{p}$ above $p$, assume $\Sigma_{\infty/\mathfrak{q},1}=\{\tau_1,\tau_2=\sigma^{-n_{\tau_1}}\tau_1,\ldots,\tau_k=\sigma^{-n_{\tau_{k-1}}}\tau_{k-1}\}$ with $\tau_1=\sigma^{-n_{\tau_k}}\tau_k$. According to Lemma \ref{calculation}, in Pic($X_T)_{\mathbb{Q}}$,
\begin{equation}
\begin{aligned}
h_{\tau_1}\  \text{vanishes on $X_T$}\quad &\Rightarrow \quad p^{n_{\tau_1}}[\omega_{\tau_2}]+[\omega_{\tau_1}]=0;\\
h_{\tau_2}\  \text{vanishes on $X_T$}\quad &\Rightarrow \quad p^{n_{\tau_2}}[\omega_{\tau_3}]+[\omega_{\tau_2}]=0;\\
&\cdots\\
h_{\tau_k}\  \text{vanishes on $X_T$}\quad &\Rightarrow \quad p^{n_{\tau_k}}[\omega_{\tau_1}]+[\omega_{\tau_k}]=0.
\end{aligned}
\end{equation}
Let $K$ be a $k\times k$ matrix with $(i,j)$-th entry given by
\begin{equation}
	K_{i,j}=\left\{
	\begin{aligned}
		1\ \  ,\quad\  &\textnormal{if } i=j,\\
		p^{n_{\tau_j}},\ \quad &\textnormal{if } i=j+1,\\
		0\ \  ,\quad\  &\textnormal{otherwise}.
	\end{aligned}\right.
\end{equation}
It is direct to check that $K$ is invertible, so $[\omega_{\tau_i}]=0$ in Pic($X_T)_{\mathbb{Q}}$ for all $1\le i\le k$.

Denote by $j:X_T\hookrightarrow X$ the embedding. By \citep{Lan} Theorem 7.2.4.1, the Hodge bundle det($\omega$) is ample on $X$, hence on $X_T$. In $\textnormal{Pic}(X_T)_{\mathbb{Q}}$ we have
\begin{equation}
\begin{aligned}
j^\ast[\text{det}(\omega)]&=[\bigotimes_{\tau'\in\Sigma_{\infty}}(\text{det}(\omega_{\mathcal{A}^\vee/X,\tilde\tau'})\otimes\text{det}(\omega_{\mathcal{A}/X,\tilde\tau'^c}))]\\
&=2\sum_{\tau'\in\Sigma_{\infty,1}}[\omega_{\mathcal{A}/X,\tilde\tau'}\otimes \omega_{\mathcal{A}/X,\tilde\tau'^c}]=2\sum_{\tau'\in\Sigma_{\infty,1}}[\omega_\tau']=2[\omega_{\tau}].
\end{aligned}
\end{equation}
Since we also know that $j^\ast[\omega^{\underline{t}}]=t_{\tau}[\omega_{\tau}]$, it follows that $t_{\tau}>0$.

\textbf{Case 2}: $\Sigma_{\infty/\mathfrak{p},1}\nsupseteq\{\tau\}$. Consider $T=\{\tau\}$, then $T'=\{\tau,\sigma^{-n_\tau}\tau\}$, $I_T=\{\sigma^{-n_\tau}\tau\}$ and $X_T$ is isomorphic to the $\mathbb{P}^1$-bundle $\mathbb{P}(\text{H}_1^{\text{dR}}(\mathcal{A}'/X')_{\sigma^{-n_\tau}\tilde\tau})$ over $X'$ by Theorem \ref{description}, with $X'=Sh_{K'}(G_{\mathcal{P}'})_k$ for the new partition $\mathcal{P}'$ and $\mathcal{A}'/X'$ the universal abelian variety. Let $\pi:X\to X'$ be the projection. Fix a closed point $s\in X'$ and write $\pi^{-1}(s)=\mathbb{P}^1_s$. Since $\omega_{\mathcal{A}^\vee/X,\sigma^{-n_\tau}\tilde\tau}$ is the tautological subbundle, we have
\begin{equation}
\omega_{\mathcal{A}^\vee/X,\sigma^{-n_\tau}\tilde\tau}|_{\mathbb{P}^1_s}\simeq\mathcal{O}_{\mathbb{P}_s^1}(-1).
\end{equation}
According to Lemma \ref{calculation},
\begin{equation}
h_\tau\ \text{vanishes on }\mathbb{P}^1_s\quad \Rightarrow \quad p^{n_\tau}[\omega_{\sigma^{-n_\tau}\tau}]|_{\mathbb{P}^1_s}+[\omega_{\tau}]|_{\mathbb{P}^1_s}=0\quad \Rightarrow\quad \omega_{\mathcal{A}^\vee/X,\tilde\tau}|_{\mathbb{P}_s^1}\simeq \mathcal{O}_{\mathbb{P}_s^1}(p^{n_\tau}) .
\end{equation}
For $\tau'\notin T'$, $\omega_{\mathcal{A}^\vee/X,\tilde\tau'}\simeq \pi^\ast\omega_{\mathcal{A}'/X',\tilde\tau'}$ by the construction in the proof of Theorem \ref{description}, so $\omega_{\mathcal{A}^\vee/X,\tilde\tau'}|_{\mathbb{P}_s^1}\simeq\mathcal{O}_{\mathbb{P}_s^1}$. Therefore,
\begin{equation}
[\omega^{\underline{t}}]|_{\mathbb{P}_s^1}=[\mathcal{O}_{\mathbb{P}_s^1}(p^{n_\tau}t_\tau-t_{\sigma^{-n_\tau}\tau})]
\end{equation}
Ampleness of $[\omega^{\underline{t}}]$ implies that $p^{n_\tau}t_\tau>t_{\sigma^{-n_\tau}\tau}$.

Thus we have completed the proof of necessity part of Theorem \ref{Main Theorem}.

\section{Sufficiency of the ampleness criterion}\label{Sufficiency}

This section is devoted to proving the sufficiency part of Theorem \ref{Main Theorem}. Our proof has three main parts: First, we reduce the proof of ampleness of $[\omega^{\underline{t}}]$ to the proof of nefness of its $\mathfrak{p}$-component $[\omega^{\underline{t}}]_\mathfrak{p}$ to be defined later; Second, we prove the non-negativity of the intersection numbers $(C\cdot [\omega^{\underline{t}}]_{\mathfrak{p}})$ for $C$ not contained in any GO-stratum; Third, for $C$ lying in some GO-stratum, we apply the description of the stratum in the previous section and use induction to reduce to the case of lower dimensional Shimura varieties.

\begin{proof} The proof is divided into five steps.

\medskip
\noindent
\textbf{Step 1} (Reduction) Recall the nef version Theorem \ref{variant}.

\begin{lemma}
If Theorem \ref{variant} holds, then the sufficiecy part of Theorem \ref{main} holds.
\begin{proof}
Let $\underline{t}$ be a tuple satisfying inequality (\ref{inequality}), i.e., $p^{n_\tau}t_\tau>t_{\sigma^{-n_\tau}\tau}$ for all $\tau\in\Sigma_{\infty,1}$. According to \citep{Lan} Theorem 7.2.4.1, the Hodge line bundle $\text{det}(\omega)$ is ample. By lemma \ref{Picard class},
\begin{equation}
\begin{aligned}
[\text{det}(\omega)] &=[\bigotimes_{\tau\in\Sigma_\infty}(\text{det}(\omega_{\mathcal{A}^\vee/X,\tilde\tau}) \otimes\text{det}(\omega_{\mathcal{A}^\vee/X,\tilde\tau^c}))]\\
&=[\bigotimes_{\tau\in\Sigma_{\infty,1}}(\omega_{\mathcal{A}^\vee/X,\tilde\tau}\otimes \omega_{\mathcal{A}^\vee/X,\tilde\tau^c})]\\
&=2\sum_{\tau\in\Sigma_{\infty,1}}[\omega_{\tau}].
\end{aligned}
\end{equation}
If we put $t'_\tau=t_\tau-2\varepsilon$, for $\varepsilon>0$ sufficiently small, this new tuple $\underline{t}'=(t'_\tau)_{\tau\in\Sigma_{\infty,1}}$ satisfies $p^{n_\tau}t_\tau\ge t_{\sigma^{-n_\tau}\tau}$ for all $\tau\in\Sigma_{\infty,1}$. But
\begin{equation}
[\omega^{\underline{t}}]=\varepsilon[\text{det}(\omega)]+[\omega^{\underline{t}'}].
\end{equation}
is the sum of an ample divisor and a nef divisor. By Proposition \ref{ampleness}, $[\omega^{\underline{t}}]$ is ample.
\end{proof}
\end{lemma}

Next, We write the sum as $[\omega^{\underline{t}}]=\sum_{\mathfrak{p}|p}[\omega^{\underline{t}}]_{\mathfrak{p}}$ with $[\omega^{\underline{t}}]_{\mathfrak{p}}=\sum_{\tau\in\Sigma_{\infty/\mathfrak{p},1}}[\omega_\tau]$. It suffices to prove each $\mathfrak{p}$-part bundle $[\omega^{\underline{t}}]_{\mathfrak{p}}$ is nef under condition (\ref{inequality'}).

\medskip
\noindent
\textbf{Step 2} (Case for a general curve) We prove that $(C\cdot [\omega^{\underline{t}}]_\mathfrak{p})\ge0$ when the curve $C$ is "generic". Recall that 
\begin{equation}
h_{\tilde\tau}\in H^0(X,\omega_{A^\vee/X,\tilde\tau}^{-1}\otimes\omega^{\otimes p^{n_\tau}}_{A^\vee/X,\sigma^{-n_\tau}\tilde\tau}).
\end{equation}
So the class of $h_\tau$ in $\text{Pic}(X)_\mathbb{Q}$ is
\begin{equation}
[h_\tau]=p^{n_\tau}[\omega_{\sigma^{-n_\tau}\tau}]-[\omega_\tau].
\end{equation}
We rewrite $[\omega^{\underline{t}}]_{\mathfrak{p}}$ as a positive linear combination of partial Hasse invariants. Let $\Sigma_{\infty/\mathfrak{p},1}=\{\tau_1,\tau_2=\sigma^{-n_{\tau_1}}\tau_1,\dots,\tau_N=\sigma^{-n_{\tau_{N-1}}}\tau_{N-1}\}$ be the ordered set with $\sigma^{-n_{\tau_N}}\tau_N=\tau_{N+1}$. Here and later, all subscripts will be viewed modulo $N$. By (\ref{Hasse}), we may write
\begin{equation}
[\omega^{\underline{t}}]_\mathfrak{p}=\sum_{i=1}^N t_i[\omega_{\tau_i}]= \sum_{i=1}^N\lambda_i [h_{\tau_i}].
\end{equation}
for $\lambda_1,\dots,\lambda_N\in\mathbb{Q}$ satisfying the following relations
\begin{equation}
\left(
\begin{matrix}
-1 & 0 & 0 & \cdots &p^{n_{\tau_N}} \\
p^{n_{\tau_1}} & -1 & 0 & \cdots &0 \\
0 & p^{n_{\tau_2}} & -1 & \cdots &0 \\
\vdots & \vdots &\vdots & \ddots &\vdots \\
0 & 0 & 0 & \cdots &-1
\end{matrix}\right)
\left(
\begin{matrix}
\lambda_1\\
\lambda_2\\
\lambda_3\\
\vdots\\
\lambda_N
\end{matrix}
\right)=
\left(
\begin{matrix}
t_1\\
t_2\\
t_3\\
\vdots\\
t_N
\end{matrix}
\right).
\end{equation}
Denote this matrix by $H$, a direct calculation gives
\begin{equation}
(H^{-1})_{i,j}= \frac{p^{n_{\tau_1}+n_{\tau_2}+\cdots+n_{\tau_N}}}{p^{n_{\tau_1}+n_{\tau_2}+\cdots+n_{\tau_N}}-1}
\begin{cases}
p^{-(n_{\tau_i}+\cdots+n_{\tau_{j-1}})} ,  & \text{if $i<j$}, \\
p^{-(n_{\tau_i}+\cdots+n_{\tau_{i+1}}+\cdots+n_{\tau_{N+j-1}})}, & \text{if $i\ge j$}.
\end{cases}
\end{equation}
In particular, all coordinates of $H^{-1}$ are positive, so all $\lambda_i$'s are non-negative.

Now for any complete curve $C\subseteq X$, if $C$ is not contained in the zero locus of any $h_{\tau_i}$, i.e., $X_{\tau_i}$, then the intersection number is always non-negative:
\begin{equation}
(C\cdot [\omega^{\underline{t}}]_\mathfrak{p})=(C\cdot\sum_{i=1}^N\lambda_i[h_{\tau_i}])=\sum_{i=1}^N\lambda_i(C\cdot[h_{\tau_i}])\ge0.
\end{equation}

\medskip
\noindent
\textbf{Step 3} (Induction base) We have to compute the intersection numbers $(C\cdot [\omega^{\underline{t}}]_\mathfrak{p})$ for $C\subseteq X_T$ for some $T\subseteq \Sigma_{\infty/\mathfrak{p},1}$. We make induction on the cardinality of $\Sigma_{\infty/\mathfrak{p},1}$. We assume that if $Y$ is the special fiber of a unitary Shimura variety with signature condition given by $\mathcal{P}'$ such that $\#\Sigma'_{\infty,1}<\#\Sigma_{\infty,1}$, and the tuple $(t'_{\tau})_{\tau\in\Sigma'_{\infty,1}}$ satisfies $p^{n'_\tau}t_\tau\ge t_{\sigma^{-n'_\tau}\tau}$ for all $\tau\in\Sigma'_{\infty,1}$, then $[\omega^{\underline{t'}}]$ is a nef line bundle on $Y$.

We check the induction base in this step.

If $\Sigma_{\infty/\mathfrak{p},1}=\emptyset$, there is nothing to do. 

If $\Sigma_{\infty/\mathfrak{p},1}=T=\{\tau\}$, the inequality condition (\ref{inequality'}) becomes $t_\tau\ge0$. We need to show that $[\omega^{\underline{t}}]_\mathfrak{p}=t_\tau[\omega_\tau]$ is nef. The partial Hasse invariant $[h_\tau]=(p^{n_\tau}-1)[\omega_\tau]$. Let $i:X_\tau\hookrightarrow X$ be the natural embedding. By lemma \ref{calculation},
\begin{equation}
h_\tau \text{ vanishes on $X_\tau$ }\Rightarrow\  (p^{n_\tau}+1)[\omega_\tau]|_{X_\tau}=0\ \Rightarrow\  [\omega_\tau]|_{X_\tau}=0.
\end{equation}
So $i^\ast[\omega^{\underline{t}}]_\mathfrak{p}=0$. The trivial class is obviously nef, hence
\begin{equation}
(C\cdot [\omega^{\underline{t}}]_{\mathfrak{p}})=(C\cdot i^\ast[\omega^{\underline{t}}]_{\mathfrak{p}})=0.
\end{equation}

\medskip
\noindent
\textbf{Step 4} (Induction step: adjacent case) Now we assume that $\#\Sigma_{\infty/\mathfrak{p},1}\ge2$.

Suppose that there exists $\{i_1,\dots,i_k\}\subseteq\Sigma_{\infty/\mathfrak{p},1}$ (in the notation of Step 2) such that $C$ is contained in the common zero locus of the partial Hasse invariants $\{h_{\tau_{i_a}}\}_{1\le a\le k}$ but not in the zero locus of any other partial Hasse invariant at $\mathfrak{p}$. We distinguish two cases: There exist adjacent labels $|i_a-i_b|\equiv1 \textnormal{ mod $N$}$ for some $a,b$ (called the adjacent case), or none of the labels are adjacent (called the sparse case). We prove in this step $(C\cdot [\omega^{\underline{t}}]_{\mathfrak{p}})\ge0$ in the adjacent case.

Without loss of generality, we may assume that $i_1=1$ and $i_2=2$.

According to lemma \ref{calculation}, in Pic($X_{\{\tau_1,\tau_2\}})_{\mathbb{Q}}$,
\begin{equation}
\begin{aligned}
h_{\tau_1}\ \textrm{vanishes on $X_{\{\tau_1,\tau_2\}}$}\quad \Rightarrow\quad [\omega_{\tau_1}]+p^{n_{\tau_1}}[\omega_{\tau_2}]=0,\\
h_{\tau_2}\ \textrm{vanishes on $X_{\{\tau_1,\tau_2\}}$}\quad \Rightarrow\quad [\omega_{\tau_2}]+p^{n_{\tau_2}}[\omega_{\tau_3}]=0.
\end{aligned}
\end{equation}
Here $[\omega_{\tau_1}]$ actually means $[\omega_{\tau_1}]|_{X_{\{\tau_1,\tau_2\}}}$, and similarly for the the other three classes. We omit this subscript from now on for simplicity. We will point out before the equation that we make the computation on which stratum.

If we denote by $j:X_{\{\tau_1,\tau_2\}}\hookrightarrow X$ the natural embedding, then
\begin{equation} 
j^\ast[\omega^{\underline{t}}]_\mathfrak{p}=(p^{n_{\tau_1}+n_{\tau_2}}t_1-p^{n_{\tau_2}}t_2+t_3)[\omega_{\tau_3}]+\sum_{i=4}^Nt_i[\omega_{\tau_i}].
\end{equation}

By Theorem \ref{description}, This GO stratum $X_{\{\tau_1,\tau_2\}}$ is isomorphic to $Y=Sh_{K'}(G_{\mathcal{P}'})_k$, the special fiber of another unitary Shimura variety. In the new partition $\mathcal{P}'$, the signature 1 part $\Sigma'_{\infty/\mathfrak{p},1}=\Sigma_{\infty/\mathfrak{p},1}\backslash \{\tau_1,\tau_2\}$. Let $\mathcal{A'}/Y$ be the universal abelian scheme. We have $\omega_{\mathcal{A}'^\vee/Y,\tilde\tau}\simeq j^\ast\omega_{\mathcal{A}/X,\tilde\tau}$ for all $\tau\in\Sigma'_{\infty/\mathfrak{p},1}$ by Remark \ref{identifications}. So $j^\ast[\omega^{\underline t}]_{\mathfrak{p}}$ is isomorphic to the automorphic line bundle $[\omega_Y^{\underline{t'}}]_{\mathfrak{p}}$ over $Y$ with $\underline t=(t_\tau)_{\tau\in\Sigma'_{\infty,1}}$. The numbers $n_\tau$'s for $Y$ is the same as $X$ except that $n'_{\tau_N}=n_{\tau_N}+n_{\tau_1}+n_{\tau_2}$. It is direct to check that
\begin{equation}
\begin{aligned}
&p^{n_{\tau_3}}(p^{n_{\tau_1}+n_{\tau_2}}t_1-p^{n_{\tau_2}}t_2+t_3)\ge p^{n_{\tau_3}}t_3\ge t_4; \\
&p^{n_{\tau_N}+n_{\tau_1}+n_{\tau_2}}t_N\ge p^{n_{\tau_1}+n_{\tau_2}}t_1\ge p^{n_{\tau_1}+n_{\tau_2}}t_1-p^{n_{\tau_2}}t_2+t_3.
\end{aligned}
\end{equation}
So $\underline{t}'$ still verifies the coefficient condition (\ref{inequality'}) for $\Sigma_{\infty,1}'$. Since $\#\Sigma'_{\infty,1}<\#\Sigma_{\infty,1}$, by the induction hypothesis, $[\omega_Y^{\underline{t}'}]$ is nef on $Y$, so the intersection number is nonnegative:
\begin{equation}
(C\cdot [\omega^{\underline{t}}]_{\mathfrak{p}})=(C\cdot j^\ast[\omega^{\underline{t}}]_{\mathfrak{p}})=(C\cdot [\omega^{\underline{t}'}]_{\mathfrak{p}})\ge0.
\end{equation}

\begin{remark} The above computation is still valid when $\#\Sigma_{\infty/\mathfrak{p},1}\le3$.
\begin{itemize}
	\item If $\Sigma_{\infty/\mathfrak{p},1}$ has cardinality $2$, then the vanishing of $h_{\tau_1}$ and $h_{\tau_2}$ forces $[\omega_{\tau_1}]=[\omega_{\tau_2}]=0$ in $\rm Pic$$(X_{\{\tau_1,\tau_2\}})_{\mathbb{Q}}$, and $(C\cdot [\omega^{\underline{t}}]_{\mathfrak{p}})=0$.
	\item If $\Sigma_{\infty/\mathfrak{p},1}$ has cardinality $3$, then $j^\ast[\omega^{\underline{t}}]_{\mathfrak{p}}=(p^{n_{\tau_1}+n_{\tau_2}}t_1-p^{n_{\tau_2}}t_2+t_3)[\omega_{\tau_3}]$. We arrive at the induction base. Since $p^{n_{\tau_1}+n_{\tau_2}}t_1-p^{n_{\tau_2}}t_2+t_3\ge0$, we have $(C\cdot [\omega^{\underline t}]_{\mathfrak{p}})\ge0$.
\end{itemize}
\end{remark}

\medskip
\noindent
\textbf{Step 5} (Induction step: sparse case) To simplify the notation, we write $\omega_i$ for $\omega_{\tau_i}$, $n_i$ for $n_{\tau_i}$ and $h_i$ for $h_{\tau_i}$. Suppose $C\subseteq X_T$ for $T=\{\tau_{i_1},\cdots,\tau_{i_k} \}$ but $C\nsubseteq X_\tau$, for any $\tau\in \Sigma_{\infty/\mathfrak{p},1}\backslash T$. We assume that $i_1<i_2<\cdots< i_k$ and $i_j-i_{j-1}\ne1$ for all $j$. This gives equations in Pic($X_T)_{\mathbb{Q}}$ by Lemma \ref{calculation}:
\begin{equation}
\begin{aligned}
h_{i_1}\ \textnormal{vanishes on $X_T$}\quad &\Rightarrow\quad [\omega_{i_1}]+p^{n_{i_1}}[\omega_{i_1+1}]=0;\\
h_{i_2}\ \textnormal{vanishes on $X_T$}\quad &\Rightarrow\quad [\omega_{i_2}]+p^{n_{i_2}}[\omega_{i_2+1}]=0;\\
&\cdots\cdots\\
h_{i_k}\ \textnormal{vanishes on $X_T$}\quad &\Rightarrow\quad [\omega_{i_k}]+p^{n_{i_k}}[\omega_{i_k+1}]=0.
\end{aligned}
\end{equation}

Denote by $i:X_T\hookrightarrow X$ the natural embedding. We have
\begin{equation}
\begin{aligned}
i^\ast[\omega^{\underline{t}}]_\mathfrak{p}&=\Big((t_{i_1}-\frac{1}{p^{n_{i_1}}}t_{i_1+1})[\omega_{i_1}]+t_{i_1+2}[\omega_{i_1+2}]+\cdots + t_{i_2-1}[\omega_{i_2-1}]\Big)\\
&+\Big((t_{i_2}-\frac{1}{p^{n_{i_2}}}t_{i_2+1})[\omega_{i_2}]+t_{i_2+2}[\omega_{i_2+2}]+\cdots + t_{i_3-1}[\omega_{i_3-1}]\Big)\\
&+\cdots\\
&+\Big((t_{i_k}-\frac{1}{p^{n_{i_k}}}t_{i_k+1})[\omega_{i_k}]+t_{i_k+2}[\omega_{i_k+2}]+\cdots + t_{i_1-1}[\omega_{i_1-1}]\Big).
\end{aligned}
\end{equation}

According to Theorem \ref{description}, for each $j=1,2,\dots,k$, the stratum $X_{\{\tau_{i_j}\}}$ is isomorphic to a $\mathbb{P}^1$-bundle over $Y_j=Sh_{K_j}(G_{\mathcal{P}_j})_k$. Denote by $\pi_j$ for this projection. In the new partition $\mathcal{P}_j$ corresponding to $Y_j$, $\Sigma^{(j)}_{\infty/\mathfrak{p},1}=\Sigma_{\infty/\mathfrak{p},1}\backslash\{\tau_{i_j},\tau_{i_j+1}\}$ and $\Sigma^{(j)}_{\infty/\mathfrak{q},1}=\Sigma_{\infty/\mathfrak{q},1}$ for all $\mathfrak{q}\ne \mathfrak{p}$. Define a tuple of numbers $(s^{(j)}_\tau)_{\tau\in \Sigma_{\infty,1}^{(j)}}$ such that $s^{(j)}_\tau=t_\tau$ for all $\tau\in\Sigma^{(j)}_{\infty,1}$. This new tuple satisfies the inequality condition for $\Sigma_{\infty,1}^{(j)}$, so by inductive hypothesis, the line bundle $[\omega^{{\underline{s}}^{(j)}}_{Y_j}]_{\mathfrak{p}}=\sum_{\tau\in\Sigma^{(j)}_{\infty/\mathfrak{p},1}}s^{(j)}_\tau[\omega_{Y_j,\tau}]$ is nef on $Y_j$. By Proposition \ref{pullback}, $\pi_j^\ast[\omega^{\underline{s}^{(j)}}_{Y_j}]_\mathfrak{p}$ is nef on $X_{\tau_{i_j}}$.

We directly compute
\begin{equation}
\begin{aligned}
\pi_j^\ast[\omega^{\underline{s}^{(j)}}_{Y_j}]_\mathfrak{p}\Big|_{X_T}&=\Big((t_{i_1}-\frac{1}{p^{n_{i_1}}}t_{i_1+1})[\omega_{i_1}]+t_{i_1+2}[\omega_{i_1+2}]+\cdots + t_{i_2-1}[\omega_{i_2-1}]\Big)\\
&+\cdots\\
&+\Big((t_{i_{j-1}}-\frac{1}{p^{n_{i_{j-1}}}}t_{i_{j-1}+1})[\omega_{i_{j-1}}]+t_{i_{j-1}+2}[\omega_{i_{j-1}+2}]+\cdots + t_{i_j-1}[\omega_{i_j-1}]\Big)\\
&+\Big(t_{i_{j}+2}[\omega_{i_{j}+2}]+\cdots + t_{i_{j+1}-1}[\omega_{i_{j+1}-1}]\Big)\\
&+\Big((t_{i_{j+1}}-\frac{1}{p^{n_{i_{j+1}}}}t_{i_{j+1}+1})[\omega_{i_{j+1}}]+t_{i_{j+1}+2}[\omega_{i_{j+1}+2}]+\cdots + t_{i_{j+2}-1}[\omega_{i_{j+2}-1}]\Big)\\
&+\cdots\\
&+\Big((t_{i_k}-\frac{1}{p^{n_{i_k}}}t_{i_k+1})[\omega_{i_k}]+t_{i_k+2}[\omega_{i_k+2}]+\cdots + t_{i_1-1}[\omega_{i_1-1}]\Big).
\end{aligned}
\end{equation}
This only differs from $i^\ast[\omega^{\underline{t}}]_\mathfrak{p}$ at the $j$-th summand.

Our strategy is to write $i^\ast[\omega^{\underline{t}}]_\mathfrak{p}$ as a positive linear combination of some partial Hasse invariants and the pullback of line bundles $[\omega^{\underline{s}^{(j)}}_{Y_j}]_\mathfrak{p}$ on $Y_j$ for each $1\le j\le k$. 

Consider the set index set $H=\{i_{1}+2,\cdots,i_2-1\}\cup\{i_2+2,\cdots,i_3-1\}\cup\cdots\cup\{i_k+2,\cdots,i_1-1\}$ and choose the Hasse invariant $h_a$ for each $a\in H$. We point out that if $i_{j+1}-i_j=2$, then the corresponding set for this part is empty. An alternative way to understand this choice is that we exactly forget those vanishing partial Hasse invariants and those partial Hasse invariants ``after'' them in the order of the $p$-adic embeddings at $\mathfrak{p}$.

Consider the equation
\begin{equation}
\begin{aligned}
i^\ast[\omega^{\underline{t}}]_\mathfrak{p}=\sum_{j=1}^k A_j\pi_j^\ast[\omega^{\underline{s}^{(j)}}_{Y_j}]_\mathfrak{p}\Big|_{X_T}+\sum_{a\in H}B_a[h_a].
\end{aligned}
\end{equation}
where the $A_j$'s and $B_a$'s are indeterminants. Expand both sides and re-group the equations properly to get a system of equations for each $j$:
\begin{equation}
\begin{aligned}
&\quad\quad t_{i_{j-1}+2}[\omega_{i_{j-1}+2}]+\cdots + t_{i_j-1}[\omega_{i_j-1}]+(t_{i_j}-\frac{1}{p^{n_{i_j}}}t_{i_j+1})[\omega_{i_j}]\\
&=\sum_{m\ne j}A_m\Big(t_{i_{j-1}+2}[\omega_{i_{j-1}+2}]+\cdots + t_{i_j-1}[\omega_{i_j-1}]+(t_{i_j}-\frac{1}{p^{n_{i_j}}}t_{i_j+1})[\omega_{i_j}]\Big)\\
&+\quad\quad A_j\Big(t_{i_{j-1}+2}[\omega_{i_{j-1}+2}]+\cdots + t_{i_j-1}[\omega_{i_j-1}]\Big)\\
&+\sum_{a=i_{j-1}+2}^{i_j-1}B_a\Big(p^{n_a}[\omega_{a+1}]-[\omega_a]\Big).
\end{aligned}
\end{equation}
Comparing the coefficients of the $[\omega]$'s, we can express the relations in matrix form.
\begin{equation}\label{matrix equation}
\left(\begin{matrix}
t_{i_{j-1}+2}\\
\vdots\\
t_{i_j-1}\\
t_{i_j}-\frac{1}{p^{n_{i_j}}}t_{i_j+1}
\end{matrix}\right)
=\sum_{m\ne j}A_m
\left(\begin{matrix}
t_{i_{j-1}+2}\\
\vdots\\
t_{i_j-1}\\
t_{i_j}-\frac{1}{p^{n_{i_j}}}t_{i_j+1}
\end{matrix}\right)
+A_j
\left(\begin{matrix}
t_{i_{j-1}+2}\\
\vdots\\
t_{i_j-1}\\
0
\end{matrix}\right)
+H_j
\left(\begin{matrix}
B_{i_{j-1}+2}\\
\vdots\\
B_{i_j-1}\\
0
\end{matrix}\right).
\end{equation}
where (the reader should be careful with the last row and column)
\begin{equation}
H_j=
\left(\begin{matrix}
-1 &0  &0  &\cdots &0 &0\\
p^{n_{i_{j-1}+2}} &-1  &0 &\cdots &0 &0\\
0 &p^{n_{i_{j-1}+3}} &-1 &\cdots &0 &0\\
\vdots &\vdots &\vdots &\ddots &\vdots &0\\
0 & 0 &0 &\cdots &-1 &0\\
0 &0 &0 &\cdots & p^{n_{i_j-1}} &0
\end{matrix}\right).
\end{equation}
We apply row operations to turn $H_j$ into its row echelon form: Add $p^{n_{i_{j-1}+2}}$ times the first row to the second row, then add $p^{n_{i_{j-1}+3}}$ times the second row to the third row, $\dots$, and finally add $p^{n_{i_j}}$ times the second last row to the last row. Comparing the last row of both sides of (\ref{matrix equation}) after these operations, we get
\begin{equation}\label{UV}
U_j+V_j=(U_j+V_j)(S-A_j)+U_jA_j.
\end{equation}
Here,
\begin{equation}
\begin{aligned}
&U_j=p^{n_{i_{j-1}+2}+\cdots+n_{{i_j}-1}}t_{i_{j-1}+2}+\cdots+p^{n_{{i_j}-1}}t_{i_j-1},\\
&V_j=t_{i_j}-\frac{1}{p^{n_{i_j}}}t_{{i_j}+1},\\
&S=\sum_{i=1}^k A_i.
\end{aligned}
\end{equation}
A special remark is that, if $i_j-i_{j-1}=2$, we do not have to choose partial Hasse invariants for this part. We only have the last row in the formula. In what follows, all we need is that $U_j\ge0$ and $V_j>0$ following from the inequality assumption at the beginning. Rearranging terms in (\ref{UV}), we get equations for all $1\le j\le k$:
\begin{equation}
\begin{aligned}
&U_jA_j=(U_j+V_j)(1+A_j-S),\\
&V_jA_j=(U_j+V_j)(S-1).
\end{aligned}
\end{equation}

The existence and explicit formula for the solutions $A_j$ can be easily deduced from the latter equation. We apply an alternative approach to estimate the values.

If $S-1=0$, then all $A_j=0$, but $1=S=\sum_{j=1}^kA_j=0$. This is a contradiction!

If $S-1<0$, then all $A_j<0$, but
\begin{equation}
0\ge \sum_{j=1}^k\frac{U_j}{U_j+V_j}A_j=\sum_{j=1}^k(1+A_j-S)=(k-1)(1-S)+1>1.
\end{equation}
This is absurd! So we must have $S-1>0$. Consequently, $A_j>0$ and $1+A_j-S\ge0$ for all $j$.

Finally we have to calculate the $B_m$'s. For this, we remove the last row in equation (\ref{matrix equation}), and move the $A_m$ terms to the left. We get
\begin{equation}
(1-S)
\left(\begin{matrix}
t_{i_{j-1}+2}\\
\vdots\\
t_{i_j-1}
\end{matrix}\right)
=H_j'
\left(\begin{matrix}
B_{i_{j-1}+2}\\
\vdots\\
B_{i_j-1}
\end{matrix}\right),
\end{equation}
where $H'_j$ is $H_j$ with the last column and row removed. An easy linear algebra shows that the inverse matrix of $H'_j$ has only non-positive entries. Using the fact that $1-S<0$, we deduce that all $B_a$'s are non-negative.

Finally we can compute intersection numbers for those $C$ essentially lying in this codimension $k$ sparse stratum:
\begin{equation}
(C\cdot [\omega^{\underline{t}}]_\mathfrak{p})=(C\cdot i^\ast[\omega^{\underline{t}}]_\mathfrak{p})= \sum_{j=1}^k A_j(C\cdot \pi_j^\ast[\omega^{\underline{s}^{(j)}}_{Y_j}]_\mathfrak{p}\Big|_{X_T})+\sum_{a\in H} B_a(C\cdot[h_a])\ge0,
\end{equation}
where the last inequality follows from the non-negativity of the coefficients, the nefness of the pullback bundle, and the assumption that $C\nsubseteq X_{\tau_a}$ for any $a\in H$. This concludes the inductive proof in the sparse case.

The sufficiency part of Theorem \ref{main} is proved.
\end{proof}

\section{Hilbert modular varieties}
In this section, we recall the construction of Hilbert modular schemes and the usual modular line bundles. The ampleness criterion for modular line bundles on Hilbert modular varieties will be an immediate corollary of the ampleness criterion for U(2)-Shimura varieties.

Let $F$ be a totally real field. Let $E=FE_0$, where $E_0=\mathbb{Q}(\sqrt{D})$ is an imaginary quadratic field whose ramification set is disjoint with that of $F$. This implies $\mathcal{O}_E=\mathcal{O}_{E_0}\otimes_{\mathbb{Z}}\mathcal{O}_F$. We assume that $p$ is unramified in $E$. Let $G$ be the algebraic group $\textnormal{Res}_{F/\mathbb{Q}}(\textnormal{GL}_{2,F})$ and $K_p=\textnormal{GL}_{2}(\mathcal{O}_F\otimes \mathbb{Z}_p)$ be its maximal compact open subgroup. Define a homomorphism $h:\mathbb{S}=\textnormal{Res}_{\mathbb{C}/\mathbb{R}}(\mathbb{G}_{m,\mathbb{C}})\to G_{\mathbb{R}}\simeq \prod_{\tau\in\Sigma_{\infty
}}\textnormal{GL}_{2,\mathbb{R}}$ by
\begin{equation}
	z=x+iy\longmapsto \left( 
	\begin{aligned}
		&\ \  x &y\\
		&-y &x
	\end{aligned}\right)
\end{equation}
We have the natural determinant map $\nu:G\to \textnormal{Res}_{F/\mathbb{Q}}(\mathbb{G}_{m,F})$. Let $G^\ast=\nu^{-1}(\mathbb{G}_{m,\mathbb{Q}})$, that is, the subset of matrices with determinant in $\mathbb{Q}$. The homomorphism $h:\mathbb{C}^\times\rightarrow G(\mathbb{R})$ factors through $G^\ast(\mathbb{R})$, so we obtain a Shimura datum $(G^\ast,[h])$. It is of PEL-type. We now give its moduli interpretation.

Let $\mathbb{A}$ and $\mathbb{A}_F$ denote the ring of adeles of $\mathbb{Q}$ and $F$, respectively. Let $W=\mathbb{A}_{F,f}^{(p)}e_1\oplus \mathbb{A}_{F,f}^{(p)}e_2$ be the standard two-dimensional $\mathbb{A}_{F,f}^{(p)}$-vector space equipped with the sympletic pairing $\psi_F:W\times W\rightarrow \mathbb{A}_f$ given by
\begin{equation}
	\psi_F(a_1e_1+a_2e_2,b_1e_1+b_2e_2)=\textnormal{Tr}_{F/\mathbb{Q}}(a_1b_2-a_2b_1)
\end{equation}
Let $K^\ast_p=\textnormal{GL}_2(\mathcal{O}_F\otimes \mathbb{Z}_p)\cap G^\ast(\mathbb{Q}_p)$. For any compact open subgroup $K^{\ast p}\subseteq G^\ast(\mathbb{A}_f^{(p)})$, we put $K^\ast=K^{\ast p}K^\ast_p\subseteq G^\ast(\mathbb{A}_f)$. We consider the functor $\textnormal{Sch}_{/\mathbb{Z}_p}^{\textnormal{loc noe}}\longrightarrow \textnormal{Sets}$ sending a locally noetherian $\mathbb{Z}_p$-scheme $S$ to the isomorphism class of tuples $(A,\iota,\lambda,\eta)$, where
\begin{itemize}
	\item $A/S$ is an abelian scheme of dimension $[F:\mathbb{Q}]$ with an $\mathcal{O}_F$-action $\iota:\mathcal{O}_F\to\textnormal{End}_S(A)$.
	\item $\lambda:A\to A^\vee$ is an $\mathcal{O}_F$-linear prime-to-$p$ quasi-polarization.
	\item Choosing a geometric point $\bar{s}_i$ for each connected component $S_i$ of $S$, $\eta$ is a collection of $\pi_1(S_i,\bar{s}_i)$-invariant $K^{\ast p}$-orbits of $\mathbb{A}_{F,f}^{(p)}$-linear isomorphisms $\eta_i:W\xrightarrow{\simeq} V^{(p)}(A_{\bar{s}_i})$ respecting the $\lambda$-Weil pairing in the following sense: There exists an isomorphism $\alpha_j\in \textnormal{Hom}(\mathbb{A}_f^{(p)},\mathbb{A}_f^{(p)}(1))$ such that the diagram 
\begin{equation}
	\begin{tikzcd}
		W\quad\ \times\quad\ W\arrow[rr, "{\langle\cdot,\cdot\rangle}"] \arrow[d, "\eta_i", shift left=10] \arrow[d, "\eta_i"', shift right=10] &  & {\mathbb{A}}_f^{(p)} \arrow[d, "\alpha_i"] \\
		V^{(p)}(A_{\bar{s}_i})\times V^{(p)}(A_{\bar{s}_i}) \arrow[rr, "\textnormal{Weil pairing}"]                                                                                    &  & {\mathbb{A}}_f^{(p)}(1)              
	\end{tikzcd}
\end{equation}
is commutative.
\end{itemize}
For $K^{\ast p}$ sufficiently small, this functor is represented by a quasi-projective smooth scheme $Sh_{K^\ast}(G^\ast)$ over $\mathbb{Z}_p$.

Now we relate this Shimura variety to the $U(2)$ Shimura varieties in this paper. Define a partition $\mathcal{P}$ by
\begin{equation}
	\mathcal{P}: \Sigma_{\infty}=\Sigma_{\infty,1}.
\end{equation}
That is, we require the Hermitian form $[\cdot,\cdot]_\tau$ in Section 2 has signature $(1,1)$ for all $\tau\in \Sigma_{\infty}$. This gives a Shimura datum $(G_{\mathcal{P}},X)$. For $K=K^pK_p\subseteq G_{\mathcal{P}}(\mathbb{A}_f)$ with $K^p$ sufficiently small, recall that we have introduced the moduli interpretation of $Sh_K(G_{\mathcal{P}})$ in the previous sections.

The group $\textnormal{GL}_2(F)$ naturally embeds into $\textnormal{GL}_2(E)$. We check that this induces an embedding $G^\ast\hookrightarrow G_{\mathcal{P}}$: Take the standard pairing $\langle\alpha,\beta\rangle=\bar{\alpha}^t\left(
\begin{aligned}
	&0 &\sqrt{D}\\
	-&\sqrt{D} &0
\end{aligned}\right)\beta$ on $V=E^{\oplus 2}$. For $\left(\begin{aligned}
&a&b\\&c&d
\end{aligned}\right)\in G^\ast$, we have
\begin{equation}
	\left(\overline{\begin{aligned}
		&a&b\\&c&d
	\end{aligned}}\right)^t
	\left(\begin{aligned}
		&0&\sqrt{D}\\-&\sqrt{D}&0
	\end{aligned}\right)
	\left(\begin{aligned}
		&a&b\\&c&d
	\end{aligned}\right)=(ad-bc)
	\left(\begin{aligned}
		&0&\sqrt{D}\\-&\sqrt{D}&0
	\end{aligned}\right).
\end{equation}
The condition for $G^\ast$ that $ad-bc\in\mathbb{Q}^\times$ implies that $\left(\begin{aligned}
	&a &b\\
	&c &d
\end{aligned}\right)\in
G_{\mathcal{P}}$.

For any sufficiently small compact open subgroup $K^{\ast p}\subseteq G^\ast(\mathbb{A}_f^{(p)})$, there exists a compact open subgroup $K^p\subseteq G(\mathbb{A}_f^{(p)})$ such that $K^{\ast p}=K^p\cap G^\ast(\mathbb{A}_f^{(p)})$. We construct a morphism $\psi:Sh_{K^\ast}(G^\ast)\longrightarrow Sh_K(G_{\mathcal{P}})$ by sending $(A,\lambda,\eta)\in Sh_{K^\ast}(G^\ast)(S)$ to $(A',\lambda',\eta')\in Sh_K(G_{\mathcal{P}})$, where
\begin{itemize}
	\item $A'=A\otimes_{\mathbb{Z}}\mathcal{O}_{E_0}$ is the Serre tensor, with the natural $\mathcal{O}_E\simeq \mathcal{O}_F\otimes_{\mathbb{Z}}\mathcal{O}_{E_0}$-action.
	\item $\lambda'=\lambda\otimes 1:A'\to A'^\vee$ is the induced prime-to-$p$ quasi-polarization of $A'$. It follows from the construction that $\lambda'$ induces complex conjugation on $\mathcal{O}_E$.
	\item Choosing a geometric point $\bar{s}_i$ for each connected component $S_i$ of $S$, $\eta'$ is a collection of $\pi_1(S_i,\bar{s}_i)$-invariant $K^p$-orbit of isomorphisms $\eta'_i:V\otimes_{\mathbb{Q}}\mathbb{A}_f^{(p)}\xrightarrow{\simeq}V^{(p)}(A'_{\bar{s}_i})$ respecting the $\lambda'$-Weil pairing. It is given by taking the $K^p$-orbit of the isomorphism
	\begin{equation}
		\alpha_{K^p}:V\otimes_{\mathbb{Q}}\mathbb{A}_f^{(p)}\xrightarrow{\simeq} W\otimes_{\mathbb{Z}}\mathcal{O}_{E_0}\xrightarrow{\alpha_{K^{\ast p}}}V^{(p)}(A_{\bar{s}_i})\otimes_{\mathbb{Z}}\mathcal{O}_{E_0}\xrightarrow{\simeq} V^{(p)}(A'_{\bar{s}_i}).
	\end{equation}
\end{itemize}

Let $X$ and $Y$ denote the special fibers of $Sh_{K}(G_{\mathcal{P}})$ and $Sh_{K^\ast}(G^\ast)$,  $\mathcal{A}$ and $\mathcal{A}^\ast$ the universal Abelian schemes, respectively. We still use $\psi$ to denote the morphism of the special fibers. We have $\mathcal{A}^\ast\otimes_{\mathbb{Z}}\mathcal{O}_{E_0}\simeq\psi^\ast\mathcal{A}$, which is non-canonically isomorphic to $\mathcal{A}^\ast\oplus\mathcal{A}^\ast$. The $\mathcal{O}_F$-action on $\mathcal{A}^\ast$ and the $\mathcal{O}_E$-action on $\mathcal{A}$ induce decompositions
\begin{equation}
	\begin{aligned}
		&\omega_{\mathcal{A}^\ast/Y}\simeq\bigoplus_{\tau\in\Sigma_{\infty}} \omega_{\mathcal{A}^\ast/Y,\tau},\\
		&\omega_{\mathcal{A}/X}\simeq\bigoplus_{\tau\in\Sigma_{\infty}} (\omega_{\mathcal{A}/X,\tilde\tau}\oplus \omega_{\mathcal{A}/X,\tilde\tau^c}).
	\end{aligned}
\end{equation}
The choice of the prime $\mathfrak{q}$ above $\mathfrak{p}$ for all $\mathfrak{p}$ over $p$ induces an identification $\omega_{\mathcal{A}^\ast/Y,\tau}\simeq \psi^\ast\omega_{\mathcal{A}/X,\tilde\tau}$ for all $\tau\in\Sigma_{\infty}$, both of which are line bundles on $Y$. This shows that $\psi$ induces an isomorphism on tangent spaces.

Let $X^{\mathrm{min}},Y^{\mathrm{min}}$ be the minimal compactifications of the special fibers. Then $\psi^{\mathrm{min}}$ is a finite $\acute{\textnormal{e}}$tale morphism. The analyse of partial Hasse invariants,  description of Goren--Oort strata, and the ampleness criterion naturally transfers to $Y^{\mathrm{min}}$:

\begin{theorem}
	Let $Y=Sh_{K^\ast}(G^\ast)_k$ be the special fiber of the Hilbert modular variety. Let $[\omega^{\underline{t}}]=\sum\limits_{\tau\in\Sigma_\infty}t_\tau[\omega_\tau]\in\textnormal{Pic}(Y)_{\mathbb{Q}}$ be the rational class associated to the tuple $\underline{t}=(t_\tau)_{\tau\in\Sigma_{\infty}}$. Then
		\begin{itemize}
			\item (Ampleness Criterion) $[\omega^{\underline{t}}]$ is an ample class if and only if $pt_\tau>t_{\sigma^{-1}\tau}$ for all $\tau\in\Sigma_{\infty}$. In particular, $t_\tau>0$ for all $\tau$.
			\item (Nefness Criterion) $[\omega^{\underline{t}}]$ is a nef class if and only if $pt_\tau\ge t_{\sigma^{-1}\tau}$ for all $\tau\in\Sigma_{\infty}$. In particular, $t_\tau\ge0$ for all $\tau$.
		\end{itemize}
\end{theorem}

\bibliographystyle{plain}
\bibliography{Ampleness}

\end{document}